\documentstyle[12pt]{article}
\newtheorem{Thm}{\indent Theorem}[section]
\newtheorem{Prop}[Thm]{\indent Proposition}
\newtheorem{Lem}[Thm]{\indent Lemma}
\newtheorem{Cor}[Thm]{\indent Corollary}

\newcommand{\Proof}{{\sl Proof.}\quad}
\newcommand{\QED}{{\unskip\nobreak\hfil\penalty50\quad\null\nobreak\hfil
{\sl q.e.d.}\parfillskip0pt\finalhyphendemerits0\par\medskip}}

\renewcommand{\H}{\mathop{\rm H}\nolimits}

\newcommand{\GM}{\mathop{\rm GM}\nolimits}
\newcommand{\CGM}{\mathop{\rm CGM}\nolimits}
\newcommand{\GEM}{\mathop{\rm GEM}\nolimits}

\newcommand{\CGEM}{\mathop{\rm CGEM}\nolimits}

\newcommand{\ic}{\mathop{\rm ic}\nolimits}
\renewcommand{\P}{\mathop{\rm P}\nolimits}
\newcommand{\Sd}{\mathop{\rm Sd}\nolimits}
\newcommand{\SdP}{\mathop{\rm SdP}\nolimits}

\newcommand{\grt}{\mathop{\rm gt}\nolimits}
\newcommand{\relint}{\mathop{\rm rel.\,int}\nolimits}

\newcommand{\Ker}{\mathop{\rm Ker}\nolimits}

\renewcommand{\Im}{\mathop{\rm Im}\nolimits}

\newcommand{\scl}{\mathrel{;}}
\newcommand{\cl}{\mathrel{:}}

\newcommand{\setmin}{{\setminus}}
\newcommand{\Zero}{{\bf 0}}

\newcommand{\calG}{{\cal G}}

\newcommand{\R}{{\bf R}}
\newcommand{\Z}{{\bf Z}}

\renewcommand{\d}{{\rm D}}
\renewcommand{\i}{{\rm i}}
\newcommand{\D}{{\bf D}}
\newcommand{\Q}{{\bf Q}}

\newcommand{\m}{{\bf m}}
\newcommand{\p}{{\bf p}}
\newcommand{\q}{{\bf q}}
\renewcommand{\t}{{\bf t}}
\renewcommand{\b}{{\bf b}}

\newcommand{\lan}{\langle}
\newcommand{\ran}{\rangle}
\newcommand{\lra}{\longrightarrow}
\newcommand{\lraw}[1]{\mathop{\lra}\limits^{#1}}
\newcommand{\ra}{\rightarrow}

\title
{Torus embeddings and\\ algebraic intersection complexes, II}

\author
{Masa-Nori ISHIDA\thanks{Supported in part by a Grant under
The Monbusho International Scientific Research Program: 04044081}
\\
{\normalsize Mathematical Institute, Faculty of Science}\\
{\normalsize Tohoku University, Sendai 980, Japan}}

\begin{document}

\maketitle%


\section*{Introduction}
\setcounter{equation}{0}

The theory of toric varieties is based on the fact that
each toric variety of dimension $r$ has an associated fan,
i.e., a finite set of rational polyhedral cones in the
$r$-dimensional real space.

The theory of intersection homologies was introduced by Goresky
and MacPherson in \cite{GM1} and \cite{GM2}.
The intersection homologies are obtained by special complexes
of sheaves on the variety, which are called the intersection
complexes.

The intersection complex of a given variety has a variation
depend on a sequence of integers which is called a perversity
\cite[\S2]{GM2}.
It is known that the complex with the middle perversity is most
important for normal complex varieties.
The decomposition theorem and strong Lefschetz theorem
for the intersection homologies with the middle perversity
were proved in \cite{BBD}.

Since these theories applied for toric varieties are used in
the combinatorial problems (cf.\cite{Stanley2}), it was
expected that these theory restricted to toric varieties are
described and proved by an combinatorial method in terms of the
associated fans (cf.\cite{Oda3}).

In \cite{Ishida3}, we introduced the additive category
$\GEM(\Delta)$ of graded exterior modules on a finite fan
$\Delta$.
We defined and constructed the intersection complex
$\ic_\p(\Delta)^\bullet$ as a finite complex in this category
for each perversity $\p$.
We defined a natural functor from the category of finite
complexes in $\GEM(\Delta)$ to that of complexes of sheaves
on the toric variety associated to the fan.
It was shown that the intersection complex of the toric
variety is obtained by applying this functor to the
intersection complex on the fan with the corresponding
perversity.
The intersection homologies of the toric varieties
are also calculated by the complex $\Gamma(\ic_\p(\Delta))^\bullet$,
where $\Gamma$ is an additive functor from $\GEM(\Delta)$
to an abelian category of graded modules over an exterior
algebra.

In this article, we work only on the fan and the category
$\CGEM(\Delta)$ of finite complexes in $\GEM(\Delta)$.

By considering the barycentric subdivision algebra of
a finte fan $\Delta$, we define an object $\SdP(\Delta)^\bullet$
of $\CGEM(\Delta)$, and we show that the intersection complex of
the fan is obtained as a quotient of $\SdP(\Delta)^\bullet$
for every perversity.

For a barycentric subdivision of a finite fan, we prove in
\S2 a decomposition theorem of the intersection complex with
the middle perversity.

In \S3, we study the case of simplicial fans.
We show that the intersection complexes of a simplicial fan
defined for any perversities between the top and the bottom
perversities are mutually quasi-isomorphic.
Some known vanishing theorems on the cohomologies of toric
varieties associated to a simplicial fan are proved in terms
of the category $\CGEM(\Delta)$.

In \S4, we get the first and the second diagonal theorems as a
consequence of the decomposition theorem in \S2
and the results for simplicial fans in \S3.

Let $\H^i(\Gamma(\ic_\p(\Delta))^\bullet)_j$ be the
homogeneous degree $j$ part of the $i$-th cohomology of
$\Gamma(\ic_\p(\Delta))^\bullet$.
Then this is a finite dimensional $\Q$-vector space
and $\H^i(\Gamma(\ic_\p(\Delta))^\bullet)_j =\{0\}$ for
$(i, j)\not\in [0, r]\times[-r, 0]$.
The first diagonal theorem says that, if $\Delta$ is a
complete fan and $\p$ is the middle perversity, then
this vector space is zero unless $i = j + r$.

The second diagonal theorem is for a cone of dimension $r$.
The strong Lefschets theorem used by Stanley \cite{Stanley1}
in the proof of $g$-conjecture on the number of faces of a
simplicial convex polytope can be replaced by this diagonal
theorem (cf.\cite[Cor.4.5]{Oda3}).


\section*{Notation}
\setcounter{equation}{0}

In the first part \cite{Ishida3} of this series of papers, the
coboundary maps of complexes are denoted by $d$ or $\partial$, and
called ``$d$-complexes'' and ``$\partial$-complexes'', respectively.
Where $\partial$ was used to represent the exterior derivatives
of logarithmic de Rham complexes and their direct sums.

In this paper, we need not use $\partial$-complexes
since we avoid to work on varieties.
Hence the $d$-complexes are simply called complexes in
this article.

The notation $E^\bullet$ means that $E$ is a complex and
its component of degree $i$ is $E^i$.
When we apply to it a functor $F$ to a category of complexes,
we usually denote by $F(E)^\bullet$ the resulting complex.
When $E^\bullet$ is a complex in an abelian category, its
$i$-the cohomology is denoted by $\H^i(E^\bullet)$.

Neither derived categories nor derived functors appear in
this article.


\section{Barycentric subdivision algebras and resolutions}
\setcounter{equation}{0}

Let $N$ be a free $\Z$-module of a fixed finite rank $r\geq 0$.
By {\em cones} in $N_\R := N\otimes_\Z\R$, we always
mean strongly convex rational polyhedral cones
(cf.\cite[Chap.1,1.1]{Oda1}).
We denote by $\Zero$ the trivial cone $\{0\}$.
For a cone $\sigma$, we set $r_\sigma :=\dim\sigma$
and $N(\sigma) := N\cap(\sigma +(-\sigma))$.
Since $\sigma$ is a rational cone, $N(\sigma)$ is a
free $\Z$-module of rank $r_\sigma$.

For each cone $\sigma$ in $N_\R$, we set
\begin{equation}
\det(\sigma) :=\bigwedge^{r_\sigma}N(\sigma)\simeq\Z\;.
\end{equation}
For cones $\sigma,\tau$ in $N_\R$ with $\sigma\prec\tau$
and $r_\tau - r_\sigma = 1$, we define the {\em incidence
isomorphism}
$q'_{\sigma/\tau}\cl\det(\sigma)\ra\det(\tau)$
as follows.

By the short exact sequence
\begin{equation}
0\lra N(\sigma)\lra N(\tau)\lra N(\tau)/N(\sigma)\lra 0
\end{equation}
of free $\Z$-modules, we get an isomorphism
\begin{equation}
\label{zzzq1}
(N(\tau)/N(\sigma))\otimes\det(\sigma)\simeq\det(\tau)\;.
\end{equation}
Since $\sigma$ is a codimension one face of $\tau$,
the image of $\tau$ in $N(\tau)_\R/N(\sigma)_\R$ is a
half line with the edge $0$.
We define the isomorphism
\begin{equation}
\label{zzzq2}
\det(\sigma)\simeq(N(\tau)/N(\sigma))\otimes\det(\sigma)
\end{equation}
by the isomorphism $\Z\ra N(\tau)/N(\sigma)$ such that $1$
is mapped into the image of $\tau$ in $N(\tau)_\R/N(\sigma)_\R$.
The isomorphism $q'_{\sigma/\tau}$ is defined to be the
composite of (\ref{zzzq2}) and (\ref{zzzq1}).
Namely, if the class $\bar a$ of $a\in N\cap\tau$ in
$N(\tau)/N(\sigma)$ is a generator, then
$q'_{\sigma/\tau}(w) = a\wedge w$ for $w\in\det(\sigma)$.

We get the following lemma (cf.\cite[Lem.1.4]{Ishida1}).


\begin{Lem}
		\label{Lem 1.1}
Let $\sigma,\rho$ be cones in $N_\R$ with $\sigma\prec\rho$
and $r_\rho - r_\sigma = 2$.
Then there exist exactly two cones $\tau$ with
$\sigma\prec\tau\prec\rho$ and $r_\tau = r_\sigma + 1$.
Let $\tau_1,\tau_2$ be the two cones.
Then the equality
\begin{equation}
q'_{\tau_1/\rho}\cdot q'_{\sigma/\tau_1} +
q'_{\tau_2/\rho}\cdot q'_{\sigma/\tau_2} = 0
\end{equation}
holds.
\end{Lem}

\Proof
The first assertion is a consequence of the fact that
any two-dimensional cone has exactly two edges
(cf.\cite[Prop.1.3]{Ishida1}).

By definiton, both $q'_{\tau_1/\rho}\cdot q'_{\sigma/\tau_1}$
and $q'_{\tau_2/\rho}\cdot q'_{\sigma/\tau_2}$ are
isomorphisms $\det(\sigma)\simeq\det(\rho)$ of the free
$\Z$-modules of rank one.
Hence it is sufficient to show that these two isomorphisms
have mutually opposite signs.
Take elements $a\in(\tau_1\setminus\sigma)\cap N$,
$b\in(\tau_2\setminus\sigma)\cap N$ and
$x\in\det(\sigma)\setminus\{0\}$.
Then $(q'_{\tau_1/\rho}\cdot q'_{\sigma/\tau_1})(x)$
has the same sign with $b\wedge a\wedge x$, while
$(q'_{\tau_2/\rho}\cdot q'_{\sigma/\tau_2})(x)$
has the same sign with $a\wedge b\wedge x$.
Hence they have distinct signs.
\QED

Let $\sigma$ and $\tau$ be nonsingular cones generated by
$\{n_1,\cdots, n_s\}$ and $\{n_1,\cdots, n_{s+1}\}$, respectively,
for a $\Z$-basis $\{n_1,\cdots, n_r\}$ of $N$ and for an integr
$0\leq s < r$ (cf.\cite[Thm.1.10]{Oda1}).
Then, it is easy to see that
\begin{equation}
q'_{\sigma/\tau}(n_1\wedge\cdots\wedge n_s) =
n_{s+1}\wedge n_1\wedge\cdots\wedge n_s =
(-1)^sn_1\wedge\cdots\wedge n_s\wedge n_{s+1}\;.
\end{equation}

If $\sigma$ and $\tau$ are simplicial cones generated
by $\{n_1,\cdots, n_s\}$ and $\{n_1,\cdots, n_{s+1}\}$,
respectively, for a $\Q$-basis $\{n_1,\cdots, n_r\}$ of $N_\Q$,
then
\begin{equation}
q'_{\sigma/\tau}(n_1\wedge\cdots\wedge n_s) =
a\cdot n_{s+1}\wedge n_1\wedge\cdots\wedge n_s
\end{equation}
for a positive rational number $a$, where we denote by
the same symbol $q'_{\sigma/\tau}$ its extension to
$\det(\sigma)\otimes\Q\ra\det(\tau)\otimes\Q$.

For a cone $\pi$ in $N_\R$, we denote by $F(\pi)$ the set
of faces of $\pi$.
The zero cone $\Zero$ and $\pi$ itself are elements of
$F(\pi)$.
For cones $\eta,\pi$ of $N_\R$ with $\eta\prec\pi$, we define
the ``closed interval''
\begin{equation}
F[\eta,\pi] :=\{\sigma\in F(\pi)\scl\eta\prec\sigma\}
\end{equation}
and the ``open interval''
\begin{equation}
F(\eta,\pi) := F[\eta,\pi]\setminus\{\eta,\pi\}\;.
\end{equation}
We also use the notation
\begin{equation}
F[\eta,\pi) := F[\eta,\pi]\setminus\{\pi\}\;.
\end{equation}

Let $\Delta$ be a finite fan \cite[1.1]{Oda1}.
We call a subset $\Phi\subset\Delta$ {\em star closed} if
$\sigma\in\Phi$ and $\sigma\prec\rho\in\Delta$ imply
$\rho\in\Phi$.
We call $\Phi$ {\em locally star closed} if
$\sigma,\rho\in\Phi$ and $\sigma\prec\rho$
imply $F[\sigma,\rho]\subset\Phi$.

Let $\Phi$ be a locally star closed subset of a finite
fan $\Delta$.
We define a finite complex $E(\Phi,\Z)^\bullet$ of free
$\Z$-modules as follows.

For each integer $i$, we set
\begin{equation}
E(\Phi,\Z)^i :=\bigoplus_{\sigma\in\Phi(i)}\det(\sigma)\;,
\end{equation}
where $\Phi(i) :=\{\sigma\in\Phi\scl r_\sigma = i\}$.
For $\sigma\in\Phi(i)$ and $\tau\in\Phi(i+1)$, the
$(\sigma,\tau)$-component of the coboundary map
\begin{equation}
d^i\cl E(\Phi,\Z)^i\lra E(\Phi,\Z)^{i+1}
\end{equation}
is defined to be $q'_{\sigma/\tau}$ if $\sigma\prec\tau$
and the zero map otherwise.
We have $d^{i+1}\cdot d^i = 0$ for every $i$ by
Lemma~\ref{Lem 1.1}.

Here we explain the relation of this complex
with the similar complex $C^\bullet(\Phi,\Z)$ in
\cite{Ishida1}.

Let $M$ be the dual $\Z$-module of $N$ and let
$M_\R := M\otimes_\Z\R$.
For a cone $\sigma$ in $N_\R$, we set
$\sigma^\perp :=\{x\in M_\R\scl\lan x, a\ran = 0,\forall a\in\sigma\}$,
where $\lan\; ,\;\ran\cl M_\R\times N_\R\ra\R$ is the natural
pairing.

For a cone $\sigma$ in $N_\R$, we set
$M[\sigma] := M\cap\sigma^\perp$, which is a free $\Z$-module
of rank $r - r_\sigma$.
We set $\Z(\sigma) :=\bigwedge^{r - r_\sigma}M[\sigma]$.

For cones $\sigma,\tau$ with $\sigma\prec\tau$ and
$r_\tau = r_\sigma + 1$, the isomorphism
\begin{equation}
q_{\sigma/\tau}\cl\Z(\sigma)\lra\Z(\tau)
\end{equation}
is defined as follows.

Let $p := r - r_\sigma$.
Then $M[\sigma]$ and $M[\tau]$ are free $\Z$-modules of
rank $p$ and $p-1$, respectively.
We take an element $n_1$ in $N$ such that the homomorphism
$\lan\;, n_1\ran\cl M[\sigma]\ra\Z$ is zero on the submodule
$M[\tau]$ and maps $M[\sigma]\cap\tau^\vee$ onto
$\{c\in\Z\scl c\geq 0\}$.
Then we define
\begin{equation}
q_{\sigma/\tau}(m_1\wedge\cdots\wedge m_p) :=
\lan m_1, n_1\ran(m_2\wedge\cdots\wedge m_p)
\end{equation}
for $m_1\in M[\sigma]$ and $m_2,\cdots, m_p\in M[\tau]$.
This definition does not depend on the choice of $n_1$.

$C^i(\Phi,\Z)$ was defined to be the direct sum
$\bigoplus_{\sigma\in\Phi(i)}\Z(\sigma)$ and the
$(\sigma,\tau)$-component of
$d^i\cl C^i(\Phi,\Z)\ra C^{i+1}(\Phi,\Z)$ was defined
to be $q_{\sigma/\tau}$ for each $(\sigma,\tau)$ with
$\sigma\prec\tau$.

For each cone $\sigma$, we define an identification
$\det(\sigma) =\Z(\sigma)\otimes\bigwedge^r N$ as follows.

Let $\{n_1,\cdots, n_r\}$ and $\{m_1,\cdots, m_r\}$ be
mutually dual $\Z$-basis of $N$ and $M$ such that
$M[\sigma]$ is generated by $\{m_1,\cdots, m_p\}$
and $N(\sigma)$ is generated by $\{n_{p+1},\cdots, n_r\}$.
Then we identify the generator
$(m_p\wedge\cdots\wedge m_1)\otimes(n_1\wedge\cdots\wedge n_r)$
of $\Z(\sigma)\otimes\bigwedge^rN$ with the generator
$n_{p+1}\wedge\cdots\wedge n_r$ of $\det(\sigma)$.

It is easy to see that $q'_{\sigma/\tau}$ and  $q_{\sigma/\tau}$
are compatible with respect to this identification.
In particular, we have
$E(\Phi,\Z)^\bullet = C^\bullet(\Phi,\Z)\otimes\bigwedge^rN$.

By this observation, we interpret some elementary
vanishing lemmas in \cite{Ishida1} and \cite{Ishida2}
to the following lemma.

Topologically, these are due to the contractability of convex sets.


\begin{Lem}
		\label{Lem 1.2}
{\rm (1)} Let $\rho$ be a cone of $N_\R$.
If $\rho\not=\Zero$, then all the cohomologies of the
complex $E(F(\rho),\Z)^\bullet$ are zero.

{\rm (2)} Let $\Delta$ be a finite fan such that the
support $\rho :=|\Delta|$ is a cone.
Let $\eta$ be an element of $\Delta$ and let
$\Phi :=\{\sigma\in\Delta\scl
\eta\prec\sigma,\sigma\cap\relint\rho\not=\emptyset\}$,
where $\relint\rho$ is the interior of $\rho$ in $N(\rho)_\R$.
Then
\begin{equation}
\H^i(E(\Phi,\Z)^\bullet) =\left\{
\begin{array}{ll}
 \{0\} & \hbox{ for } i\not= r_\rho \\
 \det(\rho) & \hbox{ for } i = r_\rho
\end{array}
\right.\;.
\end{equation}
\end{Lem}

\Proof
(1) is a special case of \cite[Prop.2.3]{Ishida1}
as it is mensioned in the article after the proposition.

For (2), let $N' := N/N(\eta)$ and let $C$ be the image of
$\rho$ in $N'_\R$.
Then $C$ is a rational polyhedral cone which is
not necessary strongly convex.
Set $\Delta(\eta{\prec}) :=\{\sigma\scl\sigma\in\Delta,\eta\prec\sigma\}$.
For each $\sigma\in\Delta(\eta{\prec})$, let $\sigma[\eta]$
be the image of $\sigma$ in $N'_\R$.
Then $\Delta[\eta] :=\{\sigma[\eta]\scl\sigma\in\Delta(\eta{\prec})\}$
is a fan of $N'_\R$ with the support $C$.
We set $\Psi :=\{\sigma[\eta]\in\Delta[\eta]\scl
\sigma[\eta]\cap\relint C\not=\emptyset\}$.
Since
$\sigma[\eta]\cap\relint C\not=\emptyset$ if and only if
$\sigma\in\Phi$, the complex $E(\Phi,\Z)^\bullet$ is naturally
isomorphic to $E(\Psi,\Z)[-r_\eta]^\bullet\otimes\det(\eta)$.
Then the lemma is a consequence of \cite[Lem.1.6]{Ishida2}
applied for $\pi = C$ and $\Sigma =\Delta[\eta]$.
Note that $C^i(\Phi,\Z_{1,0})$ in \cite{Ishida2} is equal to
$C^{i+r}(\Phi,\Z)$ for every $i$ and the coboundary
map is equal to that of $C^\bullet(\Phi,\Z)$.
Hence $\H^i(C^\bullet(\Phi,\Z_{1,0}))$ is equal to
$\H^{i+r}(C^\bullet(\Phi,\Z))$.
(In the statement of \cite[Lem.1.6]{Ishida2},
``$\pi\subset|\Sigma|\pi+(-\pi)$'' is a misprint of
``$\pi\subset|\Sigma|\subset\pi+(-\pi)$''.)
\QED

As in \cite[\S1]{Ishida3}, we denote by $A$ the exterior algebra
$\bigwedge^\bullet N_\Q$ and define the grading of $A$ by
$A_i :=\bigwedge^{-i}N_\Q$ for $i\in\Z$, where the indexes are
written as subscripts.
For a cone $\rho$, the subalgebra $\bigwedge^\bullet N(\rho)_\Q$
of $A$ is denoted by $A(\rho)$.

For a graded $\Q$-subalgebra $C\subset A$, we denote by
$\GM(C)$ the abelian category of finitely generated
left $A$-modules.
We denote by $\CGM(C)$ the abelian category of finite
complexes in $\GM(C)$.
We do not take the quotient by homotopy equivalences in the
definition of morphisms in $\CGM(C)$ in order to keep the
explicitness of the theory.
We write the complex degree by a superscript and the
graded module degree by a subscript.

Let $\Delta$ be a finite fan as before.

We consider the additive categories $\GEM(\Delta)$ and
$\CGEM(\Delta)$ as in \cite{Ishida3}.
Namely, an object $L^\bullet\in\CGEM(\Delta)$ has the following
structure.

(1) $L$ is a finite dimensional $\Q$-vector space with the
decomposition
\begin{equation}
L =\bigoplus_{\sigma\in\Delta}\;\bigoplus_{i,j\in\Z}%
L(\sigma)_j^i\;.
\end{equation}

(2) For each $\sigma\in\Delta$ and for each $i\in\Z$,
$L(\sigma)^i :=\bigoplus_{j\in\Z}L(\sigma)_j^i$ is in the
category $\GM(A(\sigma))$.

(3) For each $\sigma\in\Delta$,
$L(\sigma)^\bullet :=\bigoplus_{i\in\Z}L(\sigma)^i$ is in the
category $\CGM(A(\sigma))$.
The coboundary map is denoted by $d_L(\sigma/\sigma)$.

(4) For each pair $(\sigma,\tau)$ of distinct cones in $\Delta$
with $\sigma\prec\tau$ and for each $i\in\Z$, a homomorphism
$d_L^i(\sigma/\tau)\cl L(\sigma)^i\ra L(\tau)[1]^i = L(\tau)^{i+1}$
in $\GM(A(\sigma))$ is given, where we consider
$L(\tau)[1]^i\in\GM(A(\sigma))$ by the inclusion
$A(\sigma)\subset A(\tau)$.

(5) For each pair $(\sigma,\rho)$ of cones in $\Delta$ with
$\sigma\prec\rho$ and for each integer $i$, the equality
\begin{equation}
\label{eq in (5)}
\sum_{\tau\in F[\sigma,\rho]}
d_L^{i+1}(\tau/\rho)\cdot d_L^i(\sigma/\tau) = 0
\end{equation}
holds.

We call $\dim_\Q L$ the {\em total dimension} of $L^\bullet$.

Note that the homomoprhism
$d_L(\sigma/\tau)\cl L(\sigma)^\bullet\ra L(\tau)[1]^\bullet$
in (4) is not a homomorhism of complexes in general, i.e.,
it may not commute with the coboundary maps.
However, the condition (5) implies that this is a homomorphism
of complexes if $r_\tau - r_\sigma = 1$, since then
$F[\sigma,\rho] =\{\sigma,\rho\}$.

A homomorphism $f\cl L^\bullet\ra K^\bullet$ in $\CGEM(\Delta)$
consists of $\{f^i\cl L^i\ra K^i\scl i\in\Z\}$ which is
compatible with $d_L$ and $d_K$.
A homomorphism $f$ is said to be {\em unmixed} if
$f(\sigma/\tau)\cl L(\sigma)^\bullet\ra K(\tau)^\bullet$
is a zero map whenever $\sigma\not=\tau$.
When $f$ is unmixed, the homomorphism of complexes
$f(\sigma/\sigma)\cl L(\sigma)^\bullet\ra K(\sigma)^\bullet$
is denoted simply by $f(\sigma)$ for each $\sigma\in\Delta$.
An ummixed homomorphism has the kernel and the cokernel
in $\CGEM(\Delta)$.

For a finite fan $\Delta$, an objects
$\P(\Delta)^\bullet\in\CGEM(\Delta)$ is defined as follows.

For each $\sigma\in\Delta$, we define
$\P(\Delta)(\sigma)^\bullet :=
(\det(\sigma)\otimes A(\sigma))[-r_\sigma]^\bullet$, i.e.,
\begin{equation}
\P(\Delta)(\sigma)^i :=\left\{
\begin{array}{ll}
\det(\sigma)\otimes A(\sigma) & \hbox{ if } i = r_\sigma \\
\{0\} & \hbox{ if } i\not= r_\sigma
\end{array}
\right.\;.
\end{equation}

For $\sigma\in\Delta(i)$ and $\rho\in\Delta(i+1)$ with
$\sigma\prec\rho$, the $A(\sigma)$-homomorphism
$d^i(\sigma/\rho)\cl\P(\Delta)(\sigma)^i\ra\P(\Delta)(\rho)^{i+1}$
is defined by the homomorphism
\begin{equation}
q'_{\sigma/\rho}\otimes\lambda_{\sigma/\rho}\cl
\det(\sigma)\otimes A(\sigma)\lra \det(\rho)\otimes A(\rho)\;,
\end{equation}
where $\lambda_{\sigma/\rho}\cl A(\sigma)\ra A(\rho)$
is the natural inclusion map.
The equality (\ref{eq in (5)}) follows from
Lemma~\ref{Lem 1.1}.

As we see later, there exists an unmixed surjection
$\P(\Delta)^\bullet\ra\ic_\t(\Delta)^\bullet$ to the
intersection complex \cite[Thm.2.9]{Ishida3}
for the top perversity $\t$.
However, for a general perversity, such a description of the
intersection complex in terms of $\P(\Delta)^\bullet$ is not
possible.
We make it possible by replacing $\P(\Delta)^\bullet$ by its
{\em barycentric resolution} $\SdP(\Delta)^\bullet$.


For a finite set $\Phi$ of nontrivial cones of $N_\R$,
we define the {\em barycentric subdivision algebra} $B(\Phi)$
as follows.

We take an indeterminate $y(\sigma)$ for each $\sigma\in\Phi$,
and denote by $W(\Phi)$ the $\Q$-vector space with the basis
$\{y(\sigma)\scl\sigma\in\Phi\}$.
The graded $\Q$-algebra $B(\Phi)$ is defined to be the quotient
of the exterior $\Q$-algebra $\bigwedge^\bullet W(\Phi)$ by
the two-sided ideal generated by
\begin{equation}
\{y(\sigma)\wedge y(\rho)\scl\sigma,\rho\in\Phi\hbox{ such that }
\sigma\not\prec\rho\hbox{ and }\rho\not\prec\sigma\}\;.
\end{equation}

The grading $B(\Phi) =\bigoplus_{i=0}^\infty B(\Phi)^i$ is induced
by that of $\bigwedge^\bullet W(\Phi)$.
Note that the indexes of $B(\Phi)$ are given as superscripts.
We denote also by $y(\sigma)$ its image in $B(\Phi)^1$.
We denote by $u\cdot v$ or $uv$ the multiplicaiton
of $u, v$ in $B(\Phi)$.

We denote by $\Sd(\Phi)$ the set of sequences
$(\sigma_1,\cdots,\sigma_k)$ of distinct cones
$\sigma_1,\cdots,\sigma_k$ with
$\sigma_1\prec\cdots\prec\sigma_k$ including the case of
length zero.
For $\alpha = (\sigma_1,\cdots,\sigma_k)$ in
$\Sd(\Phi)$, we define $\max(\alpha) :=\sigma_k$
if $k > 0$ and $\max(\alpha) :=\Zero$ if $k = 0$.

For each nonnegative integer $k$, we denote by $\Sd_k(\Phi)$
the subset of $\Sd(\Phi)$ consisting of the sequences of length $k$.
By the assumption $\Zero\not\in\Phi$, we have
$\Sd_k(\Phi) =\emptyset$ for $k > r$.

For each $\alpha = (\sigma_1,\cdots,\sigma_k)$, we set
$z(\alpha) :=
y(\sigma_1)\cdots y(\sigma_k)\in B(\Phi)^k$,
where we understand $z(\alpha) = 1$ if $k = 0$.


The following lemma is clear by the definition of $B(\Phi)$.


\begin{Lem}
		\label{Lem 1.3}
For each integer $0\leq i\leq r$, the $\Q$-vector space
$B(\Phi)^i$ has the basis
\begin{equation}
\{z(\alpha)\scl\alpha\in\Sd_i(\Phi)\}\;,
\end{equation}
while $B(\Phi)^i =\{0\}$ for $i > r$.
\end{Lem}

Set $Y(\Phi) :=\sum_{\sigma\in\Phi}y(\sigma)\in B(\Phi)^1$.
By defining the coboundary map $d^i\cl B(\Phi)^i\ra B(\Phi)^{i+1}$
to be the multiplication of $Y(\Phi)$ to the left,
we regard $B(\Phi)^\bullet$ as a complex of $\Q$-vector spaces.

Let $\Delta$ be a finite fan of $N_\R$.
We consider the barycentric subdivision algebra
$B(\Delta\setmin\{\Zero\})$.
For each $\rho\in\Delta\setminus\{\Zero\}$,
$B(F(\Zero,\rho))$ is a graded subalgebra of
$B(\Delta\setmin\{\Zero\})$.

We set
$B(\rho) := B(F(\Zero,\rho))y(\rho)\subset B(\Delta\setmin\{\Zero\})$
for each $\rho\in\Delta\setmin\{\Zero\}$.
By Lemma~\ref{Lem 1.3}, $B(\rho)$ is a $\Q$-vector space
with the basis
$\{z(\alpha)y(\rho)\scl\alpha\in\Sd(F(\Zero,\rho))\}$.
We regard $B(\rho)^\bullet$ a complex by defining
$B(\rho)^i := B(F(\Zero,\rho))^{i-1}y(\rho)$ for $i\in\Z$
and defining the coboundary map to be the multiplication of
$Y(F(\Zero,\rho))$ to the left.
Note that $B(\rho)^i\not=\{0\}$ only for $1\leq i\leq r_\rho$,
where $r_\rho :=\dim\rho$.
For the zero cone $\Zero$, let $B(\Zero)^\bullet$ be the
complex defined by $B(\Zero)^0 :=\Q$ and $B(\Zero)^i :=\{0\}$
for $i\not= 0$.

For each $\beta\in\Sd(\Delta\setmin\{\Zero\})$,
$z(\beta)$ is in $B(\rho)$ if
and only if $\max(\beta) =\rho$.
If $\rho\not=\Zero$, then $\max(\beta) =\rho$ means
that $z(\beta) = z(\alpha)y(\rho)$ for some
$\alpha\in\Sd(F(\Zero,\rho))$.
Hence we get the decomposition
\begin{equation}
		\label{zzdecomp}
B(\Delta\setmin\{\Zero\}) =
\bigoplus_{\rho\in\Delta}B(\rho)
\end{equation}
as a $\Q$-vector space by Lemma~\ref{Lem 1.3}.
However this is not a direct sum of the complexes.

We introduce a decreasing filtration $\{F^k(B(\Delta\setmin\{\Zero\}))\}$ by
\begin{equation}
F^k(B(\Delta\setmin\{\Zero\})) :=
\bigoplus_{i=k}^r\bigoplus_{\rho\in\Delta(i)}B(\rho)\;.
\end{equation}
Then it is easy to see that each $F^k(B(\Delta\setmin\{\Zero\}))$ is a
subcomplex of $B(\Delta\setmin\{\Zero\})^\bullet$ for each $k$, and
\begin{equation}
F^k(B(\Delta\setmin\{\Zero\}))^\bullet/%
F^{k+1}(B(\Delta\setmin\{\Zero\}))^\bullet =
\bigoplus_{\rho\in\Delta(k)}B(\rho)^\bullet
\end{equation}
as complexes.

For $\rho,\mu\in\Delta$, with $\rho\prec\mu$ and $\rho\not=\mu$,
we define a $\Q$-linear map
\begin{equation}
\label{zzzmult}
\varphi_{\rho/\mu}\cl B(\rho)\lra B(\mu)
\end{equation}
to be the multiplication of $y(\mu)$ to the left, i.e.,
$w\mapsto y(\mu)w$.
If $r_\mu - r_\rho = 1$, then $\varphi_{\rho/\mu}$ is
a homomorphism $B(\rho)^\bullet\ra B(\mu)[1]^\bullet$
of complexes, however it is not the case if $r_\mu - r_\rho > 1$.

We consider the complex of $A(\rho)$-modules
\begin{equation}
(B(\rho)\otimes_\Q A(\rho))^\bullet :=
B(\rho)^\bullet\otimes_\Q A(\rho)\;.
\end{equation}
By definition, $B(\rho)^i\otimes_\Q A(\rho)_j\not=\{0\}$ only for
$1\leq i\leq r_\rho$ and $-r_\rho\leq j\leq 0$, if $\rho\not=\Zero$.

The barycentric resolution $\SdP(\Delta)^\bullet\in\CGEM(\Delta)$
of $\P(\Delta)^\bullet$ is defined as follows.

For each $\rho\in\Delta$, we set
\begin{equation}
\SdP(\Delta)(\rho)^\bullet :=
(B(\rho)\otimes_\Q A(\rho))^\bullet\;.
\end{equation}
In particular,
\begin{equation}
\SdP(\Delta)(\rho)^i =
\bigoplus_{\alpha\in\Sd_{i-1}(F(\Zero,\rho))}%
(\Z z(\alpha)y(\rho))\otimes A(\rho)
\end{equation}
for each $1\leq i\leq r_\rho$ if $\rho\not=\Zero$, while
$\SdP(\Delta)(\Zero)^0 =\Q$ and $\SdP(\Delta)(\Zero)^i =\{0\}$
for $i\not= 0$.

For $\rho,\mu\in\Delta$ with $\rho\prec\mu$ and $\rho\not=\mu$,
the $(\rho,\mu)$-component
\begin{equation}
d_{\SdP(\Delta)}(\rho/\mu)\cl
\SdP(\Delta)(\rho)^\bullet\lra\SdP(\Delta)(\mu)[1]^\bullet
\end{equation}
of the coboundary map of $\SdP(\Delta)^\bullet$ is defined to
be the tensor product of $\varphi_{\rho/\mu}\cl B(\rho)\ra B(\mu)$
defined at (\ref{zzzmult}) and the natural inclusion map
$A(\rho)\ra A(\mu)$.
Note that this is not a homomorphism of complexes, if
$r_\mu - r_\rho > 1$.

In order to check the equality (\ref{eq in (5)}) for
$\SdP(\Delta)^\bullet$, it is sufficient to show the
equality
\begin{equation}
\sum_{\tau\in F[\rho,\mu]}\varphi_{\tau/\mu}\cdot\varphi_{\rho/\tau}
= 0\;,
\end{equation}
where $\varphi_{\rho/\rho}$ and $\varphi_{\mu/\mu}$ are the
multiplication of $Y(F(\Zero,\rho))$ and $Y(F(\Zero,\mu))$
to the left, respectively.
This equality follows from the equality
\begin{equation}
(\;y(\mu)Y(F(\Zero,\rho)) + Y(F(\Zero,\mu))y(\mu) +
\sum_{\tau\in F(\rho,\mu)}y(\mu)y(\tau)\;)\;y(\rho) = 0
\end{equation}
in $B(\Delta\setmin\{\Zero\})$ which is checked easily.

For each $\rho\in\Delta$, we define a covariant functor
\begin{equation}
\i_\rho^\circ\cl\CGEM(\Delta)\lra\CGM(A(\rho))
\end{equation}
as follows.
For $L^\bullet\in\CGEM(\Delta)$ and for each $i\in\Z$,
we set
\begin{equation}
\i_\rho^\circ(L)^i :=
\bigoplus_{\sigma\in F[\Zero,\rho)}L(\sigma)_{A(\rho)}^i\;.
\end{equation}
For $\sigma,\tau\in F[\Zero,\rho)$ with $\sigma\prec\tau$,
the $(\sigma,\tau)$-component of
$d^i\cl\i_\rho^\circ(L)^i\ra\i_\rho^\circ(L)^{i+1}$ is defined
to be the $A(\rho)$-homomorphism induced by $d_L^i(\sigma/\tau)$.

Recall that a similar functor $\i_\rho^*$ was defined in
\cite[\S2]{Ishida3}.
We get the definition of $\i_\rho^*$ by replacing
$F[\Zero,\rho)$ in the definition of $\i_\rho^\circ$
by $F[\Zero,\rho] = F(\rho)$.

For $L^\bullet\in\CGEM(\Delta)$ and $\rho\in\Delta$, we get
an exact sequence
\begin{equation}
0\lra L(\rho)^\bullet\lra\i_\rho^*(L)^\bullet
\lra\i_\rho^\circ(L)^\bullet\lra 0
\end{equation}
in $\CGM(A(\rho))$.
Since $\i_\rho^*(L)^i = L(\rho)^i\oplus\i_\rho^\circ(L)^i$
for each integer $i$, $\i_\rho^*(L)^\bullet$ is equal to
the mapping cone of a homomoprhism
\begin{equation}
\phi(L,\rho)\cl\i_\rho^\circ(L)^\bullet\lra  L(\rho)[1]^\bullet
\end{equation}
of complexes.
The restriction of $\phi(L,\rho)^i$ to the component
$L(\sigma)_{A(\rho)}^i$ is the $A(\rho)$-homomorphism
induced by $d_L^i(\sigma/\rho)$ for each $\sigma$.

Let $f\cl L^\bullet\ra K^\bullet$ be a homomorphism
in $\CGEM(\Delta)$.
For $\rho\in\Delta$, consider the following diagram
\begin{equation}
\begin{array}{ccc}
\label{circle diagram}
\makebox[20pt]{}\i_\rho^\circ(L)^\bullet &
\mathop{\lra}\limits^{\i_\rho^\circ(f)} &
 \i_\rho^\circ(K)^\bullet\makebox[10pt]{} \\
{\phi(L,\rho)}\downarrow &
 &\downarrow{\phi(K,\rho)} \\
\makebox[20pt]{}L(\rho)[1]^\bullet &
\mathop{\lra}\limits^{f(\rho)[1]} &
K(\rho)[1]^\bullet\makebox[10pt]{}
\end{array}
..
\end{equation}
Let $\{u^i\scl i\in\Z\}$ be the collection of
$A(\rho)$-homomorphisms
$u^i\cl\i_\rho^\circ(L)^i\ra K(\rho)^i = K(\rho)[1]^{i-1}$
induced by $\{f^i(\sigma/\rho)\scl\sigma\in F[\Zero,\rho)\}$.
Then restriction of the equality $d_K\cdot f = f\cdot d_L$
to $\i_\rho^\circ(L)^i$ implies
\begin{equation}
\phi(K,\rho)^i\cdot\i_\rho^\circ(f)^i + d_{K(\rho)[1]}^i\cdot u^i =
f(\rho)[1]^i\cdot\phi(L,\rho)^i + u^{i+1}\cdot d_{\i_\rho^\circ(L)}^i\;.
\end{equation}
Hence the difference
$\phi(K,\rho)^i\cdot\i_\rho^\circ(f)^i - f(\rho)[1]^i\cdot\phi(L,\rho)^i$
is homotopy equivalent to the zero map.
In particlar, the diagram (\ref{circle diagram}) induces a
commutative diagram of cohomologies.
If $f$ is unmixed, then all $u^i$'s are zero and the
diagram (\ref{circle diagram}) is commutative.

Let $L^\bullet$ be an object of $\CGEM(\Delta)$.
$K^\bullet\in\CGEM(\Delta)$ with a homomorphism
$f\cl K^\bullet\ra L^\bullet$ is said to be a {\em subcomplex}
if $f$ is unmixed and
$f(\sigma)\cl K(\sigma)^\bullet\ra L(\sigma)^\bullet$
is an inclusion map of complexes for every $\sigma\in\Delta$.
By a {\em homogeneous element} of $L^\bullet$, we mean
an element of $L(\sigma)_j^i$ for some $\sigma\in\Delta$ and
$i, j\in\Z$.
For a set $S$ of homogeneous elements of $L^\bullet$, there
exists a unique subcomplex $\lan S\ran^\bullet$ of $L^\bullet$
generated by $S$.
The complex $\lan S\ran^\bullet$ is described inductively as follows.

For $\sigma\in\Delta$ and $i, j\in\Z$, let
$S(\sigma, i, j) := S\cap L(\sigma)_j^i$.
We denote by $\lan S(\sigma)\ran^\bullet$ the
$A(\sigma)$-subcomplex of $L(\sigma)^\bullet$ generated by
$S(\sigma) :=\bigcup_{i,j\in\Z}S(\sigma, i, j)$ for
each $\sigma\in\Delta$.

We set $\lan S\ran(\Zero)^\bullet :=\lan S(\Zero)\ran^\bullet$.
Let $\rho$ be in $\Delta\setminus\{\Zero\}$, and assume that
we already know $\lan S\ran(\sigma)^\bullet$ for
$\sigma\in F[\Zero,\rho)$.
Then $\lan S\ran(\rho)^\bullet$ is the $A(\rho)$-subcomplex
of $L(\rho)^\bullet$ given by
\begin{equation}
\label{eq in S}
\lan S\ran(\rho)^\bullet =
\phi(L,\rho)(\i_\rho^\circ(\lan S\ran))^\bullet +
\lan S(\rho)\ran^\bullet\;.
\end{equation}

Let $\p$ be a perversity of $\Delta$, i.e., $\p$ is a map
$\Delta\setminus\{\Zero\}\ra\Z$.

We denote by $k_\p(\Delta)^\bullet$ the subcomplex of
$\SdP(\Delta)^\bullet$ generated by
\begin{equation}
\bigcup_{\sigma\in\Delta\setminus\{\Zero\}}\;%
\bigcup_{i+j\leq\p(\sigma)}\SdP(\Delta)(\sigma)_j^i\;.
\end{equation}

Note that if we set this set $S$, then $\lan S(\sigma)\ran^\bullet$
is equal to the gradual truncation
$\widetilde{\grt}_{\leq\p(\sigma)}\SdP(\Delta)(\sigma)^\bullet$
(cf.\cite[\S1]{Ishida3}).


\begin{Lem}
\label{Lem 1.4}
Let $\Delta$ be a finite fan.
For each $\rho\in\Delta\setminus\{\Zero\}$, the homomorphism
\begin{equation}
\phi(\SdP(\Delta),\rho)\cl\i_\rho^\circ(\SdP(\Delta))^\bullet
\lra \SdP(\Delta)(\rho)[1]^\bullet
\end{equation}
is an isomorphism.
\end{Lem}

\Proof
Since $\SdP(\Delta)^i =\{0\}$ for $i < 0$ and
$\SdP(\Delta)(\rho)^0 =\{0\}$ for $\rho\not=\Zero$,
it is sufficient to show that $\phi(\SdP(\Delta),\rho)^i$
is an isomorphism for $i\geq 0$.
For $i = 0$, we have
$i_\rho^\circ(\SdP(\Delta))^0 =\Z\otimes A(\rho)$ and
$\SdP(\Delta)(\rho)^1 =\Z y(\rho)\otimes A(\rho)$, and
$\phi(\SdP(\Delta),\rho)^0$ is the isomorphism given by
$1\otimes 1\mapsto y(\rho)\otimes 1$.
Assume $i > 0$.
The free $A(\rho)$-module
$\SdP(\Delta)(\rho)^{i+1} = B(\rho)^{i+1}\otimes_\Q A(\rho)$,
has the basis $\{z(\alpha)y(\rho)\scl\alpha\in\Sd_i(F(\Zero,\rho))\}$
which is decomposed to the disjoint union
\begin{equation}
\bigcup_{\sigma\in F(\Zero,\rho)}
\{z(\beta)y(\sigma)y(\rho)\scl\beta\in\Sd_{i-1}(F(\Zero,\sigma))\}\;.
\end{equation}

On the other hand, the component $\SdP(\Delta)(\sigma)_{A(\rho)}^i$
of $\i_\rho^\circ(\SdP(\Delta))^i$ has the basis
\begin{equation}
\{z(\beta)y(\sigma)\scl\beta\in\Sd_{i-1}(F(\Zero,\sigma))\}
\end{equation}
for each $\sigma\in F(\Zero,\rho)$.
Since
$\phi(\SdP(\Delta),\rho)^i(z(\beta)y(\sigma)) =
(-1)^iz(\beta)y(\sigma)y(\rho)$ for each $z(\beta)y(\sigma)$,
$\phi(\SdP(\Delta),\rho)^i$ induces an isomoprhism from
$\SdP(\Delta)(\sigma)_{A(\rho)}^i$ to the submodule of
$\SdP(\Delta)(\rho)^{i+1}$ generated by
$\{z(\beta)y(\sigma)y(\rho)\scl\beta\in\Sd_{i-1}(F(\Zero,\sigma))\}$.
We are done, since $\phi(\SdP(\Delta),\rho)^i$ is the direct sum of
these isomorphisms for $\sigma\in F(\Zero,\rho)$.
\QED

We define an unmixed homomorphism
\begin{equation}
\psi_\Delta\cl\SdP(\Delta)^\bullet\lra\P(\Delta)^\bullet
\end{equation}
of complexes in $\CGEM(\Delta)$ as follows.

We define a homomorphism
\begin{equation}
\label{sdp to p}
\psi_\Delta(\rho)\cl\SdP(\Delta)(\rho)^\bullet
\lra\P(\Delta)(\rho)^\bullet
\end{equation}
in $\CGM(A(\rho))$ for each $\rho\in\Delta$.
Let $k := r_\rho$.
Then $\psi_\Delta(\rho)^i := 0$ for $i\not= k$ since
$\P(\Delta)(\rho)^i =\{0\}$.
For $\alpha\in\Sd_{k-1}(F(\Zero,\rho))$, we define
$a :=\psi_\Delta(\rho)^k(z(\alpha)y(\rho))$ to be the
generator of $\det(\rho)\subset\P(\Delta)(\rho)^k$ which
is determined by the orientation of the sequence
$\alpha = (\sigma_1,\cdots,\sigma_{k-1})$, i.e.,
if we take $a_i\in N\cap\relint\sigma_i$ for $i = 1,\cdots, k-1$
and $a_k\in N\cap\relint\rho$, then $a$ is the generator of
$\det(\rho)$ which has the same sign with
$a_1\wedge\cdots\wedge a_k\in\det(\rho)$.
The commutativity of $\psi_\Delta$ and the coboundary
maps is checked easily.
The only one nontrivial commutativity is of
the component for $\SdP(\Delta)(\rho)^{k-1}$ and
$\P(\Delta)(\mu)^k$ with $\mu =\rho$.
For each $\alpha\in\Sd_{k-2}(F(\Zero,\rho))$, there exist
exactly two $\beta_1,\beta_2\in\Sd_{k-1}(F(\Zero,\rho))$
which contains $\alpha$ as a subsequence.
We get the commutativity, since $\beta_1$ and $\beta_2$
have mutually opposite orientations and hence the equality
\begin{equation}
\psi_\Delta(\rho)^k(z(\beta_1)y(\rho)) +
\psi_\Delta(\rho)^k(z(\beta_2)y(\rho)) = 0
\end{equation}
holds.

By the definition, $\psi_\Delta(\rho)$ is surjective.
The following lemma implies that the kernel is generated by
\begin{equation}
\label{ker of sdp to p}
\bigcup_{i=1}^{r_\rho-1}\SdP(\Delta)(\rho)^i
\end{equation}
as a subcomplex for each $\rho\in\Delta\setminus\{\Zero\}$.


\begin{Lem}
		\label{Lem 1.5}
The above unmixed homomorphism $\psi_\Delta$ is a
quasi-isomorphism in $\CGEM(\Delta)$.
\end{Lem}

\Proof
Since $\SdP(\Delta)(\Zero) =\P(\Delta)(\Zero) =\Q$,
$\psi_\Delta(\Zero)$ is a quasi-isomorphism.
Let $\Phi$ be a maximal subfan of $\Delta$ such that
the restriction of $\psi_\Delta$ to $\Phi$ is a quasi-isomorphism.
Suppose $\Phi\not=\Delta$, and let $\rho$ be a minimal
element of $\Delta\setminus\Phi$.

Since $\psi_\Delta$ is unmixed, we get a commutative diagram
\begin{equation}
\begin{array}{ccc}
\makebox[70pt]{}\i_\rho^\circ(\SdP(\Delta))^\bullet &
\mathop{\lra}\limits^{\i_\rho^\circ\psi_\Delta} &
 \i_\rho^\circ(\P(\Delta))^\bullet\makebox[50pt]{} \\
{\phi(\SdP(\Delta),\rho)}\downarrow &
 &\downarrow{\phi(\P(\Delta),\rho)} \\
\makebox[70pt]{}\SdP(\Delta)(\rho)[1]^\bullet &
\mathop{\lra}\limits^{\psi_\Delta(\rho)[1]} &
 \P(\Delta)(\rho)[1]^\bullet\makebox[50pt]{}
\end{array}
..
\end{equation}

By Lemma~\ref{Lem 1.4}, $\phi(\SdP(\Delta),\rho)$ is an
isomorphism, while $\i_\rho^\circ\psi_\Delta$ is a quasi-isomorphism
since it depends only on the restriction of $\psi_\Delta$ to $\Phi$.

Since the mapping cone $\i_\rho^*(\P(\Delta))^\bullet$ of
$\phi(\P(\Delta),\rho)$ is equal to
$E(F(\rho),\Z)^\bullet\otimes A(\rho)$, it has trivial
cohomologies by Lemma~\ref{Lem 1.2},(1).
Hence $\phi(\P(\Delta),\rho)$ is also a quasi-isomorphism.

By the commutative diagram, $\psi_\Delta(\rho)$ is also a
quasi-isomorphism.
This contradicts the maximality of $\Phi$.
\QED


\begin{Lem}
\label{Lem 1.6}
For $L^\bullet\in\CGEM(\Delta)$ and a homomorphism
$f_0\cl\SdP(\Delta)(\Zero)^\bullet\ra L(\Zero)^\bullet$,
there exists a unique unmixed homomorphism
$f\cl\SdP(\Delta)^\bullet\ra L^\bullet$ with
$f(\Zero) = f_0$.
\end{Lem}

\Proof
We prove the lemma by induction on the number of
cones in $\Delta$.

If $\Delta =\{\Zero\}$, then the assertion is clear.

Assume that $\Delta\not=\{\Zero\}$ and $\pi$ is a
maximal element of $\Delta$.
Set $\Delta' :=\Delta\setminus\{\pi\}$ and assume that
$f'\cl\SdP(\Delta')^\bullet\ra (L|\Delta')^\bullet$
is the unique extension of $f_0$.
Let $f_1\cl\SdP(\Delta)(\pi)^\bullet\ra L(\pi)^\bullet$
be a homomorphism in $\CGM(A(\pi))$.
Then $f'$ is extended to an unmixed homomorphism
$f\cl\SdP(\Delta)^\bullet\ra L^\bullet$ by $f(\pi) := f_1$
if and only if the diagram
\begin{equation}
\begin{array}{ccc}
\makebox[70pt]{}\i_\pi^\circ(\SdP(\Delta))^\bullet &
\mathop{\lra}\limits^{\i_\pi^\circ f'} &
 \i_\pi^\circ(L)^\bullet\makebox[50pt]{} \\
{\phi(\SdP(\Delta),\pi)}\downarrow &
 &\downarrow{\phi(L,\pi)} \\
\makebox[70pt]{}\SdP(\Delta)(\pi)[1]^\bullet &
\mathop{\lra}\limits^{f_1[1]} & L(\pi)[1]^\bullet\makebox[50pt]{}
\end{array}
\end{equation}
is commutative.
Since $\phi(\SdP(\Delta),\pi)$ is an isomorphism by Lemma~\ref{Lem 1.4},
such a $f_1$ exists uniquely.
\QED


\begin{Thm}
\label{Thm 1.7}
Let $\p$ be a perversity of $\Delta$.
Let
\begin{equation}
\varphi(\Delta,\p)\cl
\SdP(\Delta)^\bullet\lra\ic_\p(\Delta)^\bullet
\end{equation}
be the unmixed homomorphism obtained by extending the
identity $\SdP(\Delta)(\Zero) =\ic_\p(\Delta)(\Zero) =\Q$.
Then $\varphi(\Delta,\p)$ is surjective and the kernel
is equal to $k_\p(\Delta)^\bullet$.
\end{Thm}

\Proof
We prove the theorem by induction.
Namely, let $\Phi$ be the maximal subfan of $\Delta$
such that $\varphi(\Delta,\p)(\sigma)$ is surjective
and the kernel is equal to $k_\p(\Delta)(\sigma)^\bullet$
for $\sigma\in\Phi$.
Since $k_\p(\Delta)(\Zero) =\{0\}$, we have $\Zero\in\Phi$.
Suppose $\Phi\not=\Delta$ and let $\rho$ be a minimal
element of $\Delta\setminus\Phi$.
Since $\varphi(\Delta,\p)$ is unmixed,
we get the following commutative diagram.
\begin{equation}
\begin{array}{ccccccccc}
0 & \ra & \i_\rho^\circ(k_\p(\Delta))^\bullet &
\lra & \i_\rho^\circ(\SdP(\Delta))^\bullet &
\mathop{\lra}\limits^{w} & \i_\rho^\circ(\ic_\p(\Delta))^\bullet & \ra & 0 \\
 & & \hbox{  }\downarrow{\phi_1} & & \hbox{  }\downarrow{\phi_2}
 & & \hbox{  }\downarrow{\phi_3} & & \\
0 & \ra & k_\p(\Delta)(\rho)[1]^\bullet &
\mathop{\lra}\limits^{u} & \SdP(\Delta)(\rho)[1]^\bullet &
\mathop{\lra}\limits^{v} & \ic_\p(\Delta)(\rho)[1]^\bullet & \ra & 0
\end{array}
..
\end{equation}

The upper line of this diagram is exact by the assumption.
Among the vertical homomorphisms, $\phi_2$ is an isomorhism
by Lemma~\ref{Lem 1.4} and $\phi_3$ is the natural surjection
\begin{equation}
\i_\rho^\circ(\ic_\p(\Delta))^\bullet\lra
\grt^{\geq\p(\rho)}(\i_\rho^\circ(\ic_\p(\Delta)))^\bullet
= \ic_\p(\Delta)(\rho)[1]^\bullet
\end{equation}
by the construction of $\ic_\p(\Delta)^\bullet$ \cite[Thm 2.9]{Ishida3}.
Hence $v$ is surjective, while $u$ is an inclusion map.

Since there exists an exact sequence
\begin{equation}
0\ra
\widetilde{\grt}_{\leq\p(\rho)-1}(\i_\rho^\circ(\ic_\p(\Delta)))^\bullet
\ra\i_\rho^\circ(\ic_\p(\Delta))^\bullet\ra
\grt^{\geq\p(\rho)}(\i_\rho^\circ(\ic_\p(\Delta)))^\bullet\ra 0\;,
\end{equation}
$\Ker\phi_3$ is equal to
$\widetilde{\grt}_{\leq\p(\rho)-1}(\i_\rho^\circ(\ic_\p(\Delta)))^\bullet$.
Since $w$ is surjective, this is equal to the image of
$\widetilde{\grt}_{\leq\p(\rho)-1}(\i_\rho^\circ(\SdP(\Delta)))^\bullet$
by $w$.
Hence
\begin{equation}
w^{-1}(\Ker\phi_3) =
\i_\rho^\circ(k_\p(\Delta))^\bullet +
\widetilde{\grt}_{\leq\p(\rho)-1}(\i_\rho^\circ(\SdP(\Delta)))^\bullet\;.
\end{equation}
Since
\begin{equation}
\widetilde{\grt}_{\leq\p(\rho)-1}(\SdP(\Delta)(\rho))[1])^\bullet =
(\widetilde{\grt}_{\leq\p(\rho)}(\SdP(\Delta)(\rho)))[1]^\bullet\;,
\end{equation}
$\Ker v =\phi_2(w^{-1}(\Ker\phi_3))$ is equal to the sum of $\Im\phi_1$
and $(\widetilde{\grt}_{\leq\p(\rho)}(\SdP(\Delta)(\rho)))[1]^\bullet$
in $\SdP(\Delta)(\rho)[1]^\bullet$.

On the other hand, $k_\p(\Delta)(\rho)^\bullet$ is the
subcomplex of $\SdP(\Delta)(\rho)^\bullet$ generated by
$\Im\phi_1$ and
$\bigcup_{i+j\leq\p(\rho)}\SdP(\Delta)(\rho)_j^i$ (cf.(\ref{eq in S})).
Since the subcomplex of $\SdP(\Delta)(\rho)^\bullet$
generated by the last set is
$\widetilde{\grt}_{\leq\p(\rho)}(\SdP(\Delta)(\rho))^\bullet$,
the lower line is also exact.
This contradicts the maximality of $\Phi$, and we conclude
$\Phi =\Delta$.
\QED

For a finite fan $\Delta$, the covariant additive
functor $\Gamma\cl\GEM(\Delta)\ra\GM(A)$ is defined
as follows.

For $L\in\GEM(\Delta)$, we set
\begin{equation}
\Gamma(L) :=\bigoplus_{\sigma\in\Delta}L(\sigma)_A\;.
\end{equation}
Let $f\cl L\ra K$ be a homomorphism in $\GEM(\Delta)$.
For $\sigma,\tau\in\Delta$, the $(\sigma,\tau)$-component
of the homomorphism $\Gamma(f)\cl\Gamma(L)\ra\Gamma(K)$
is defined to be the $A$-homomorphism
$L(\sigma)_A\ra K(\tau)_A$
induced by $f(\sigma/\tau)$ and is the zero map otherwise.

Let $L^\bullet$ be in $\CGEM(\Delta)$.
Then $\Gamma(L)^\bullet$ is in $\CGM(A)$, i.e.,
is a finite complex of graded $A$-modules.
For each integer $q$, the homogeneous component
$\Gamma(L)_q^\bullet$ is a complex of $\Q$-vector spaces.

Since $A(\rho)_A = A$ for all $\rho\in\Delta$, we have
\begin{equation}
\Gamma(\SdP(\Delta)^i) =
B(\Delta\setmin\{\Zero\})^i\otimes_\Q A
\end{equation}
for each $i\in\Z$.
By comparing the definitions of the coboundary maps,
we get the following
lemma.


\begin{Lem}
		\label{Lem 1.8}
For any finite fan $\Delta$,
$\Gamma(\SdP(\Delta))^\bullet$ is canonically isomorphic
to $B(\Delta\setmin\{\Zero\})^\bullet\otimes_\Q A$ as a
complex of $A$-modules.
\end{Lem}

The {\em top perversity} $\t$ is defined by $\t(\sigma) := r_\sigma - 1$
for every nontrivial cone $\sigma$, while the
{\em bottom perversity} $\b$ is defined by $\b :=-\t$.
We use this notation for all finite fans.


\begin{Lem}
		\label{Lem 1.9}
Let $\Delta$ be a finite fan.
Then the unmixed surjection $\varphi(\Delta,\t)\cl
\SdP(\Delta)^\bullet\ra\ic_\t(\Delta)^\bullet$
induces an unmixed surjection
$\varphi\cl\P(\Delta)^\bullet\ra\ic_\t(\Delta)^\bullet$ such that
\begin{equation}
\Ker\varphi(\sigma)^{r_\sigma} =
\det(\sigma)\otimes(N(\sigma)A(\sigma))
\end{equation}
for every $\sigma\in\Delta$.
\end{Lem}

Note that $N(\sigma)A(\sigma)$ is a two-sided maximal ideal
of $A(\sigma)$.

\Proof
The kernel of the homomorphism $\varphi(\Delta,\t)$ is
equal to $k_\t(\Delta)^\bullet$ by Theorem~\ref{Thm 1.7}.
By definition, this subcomplex is generated by the union of
\begin{equation}
\bigcup_{i+j<r_\sigma}\SdP(\Delta)(\sigma)_j^i
\end{equation}
for all $\sigma\in\Delta\setminus\{\Zero\}$.
Since $\SdP(\Delta)(\sigma)_j^i$ is nonzero only for
$1\leq i\leq r_\sigma,\; -r_\sigma\leq j\leq 0$,
the condition $i+j<r_\sigma$ means all
$(i, j)$ except for $(r_\sigma, 0)$.
On the other hand, the kernel of
$\psi_\Delta(\sigma)\cl
\SdP(\Delta)(\sigma)^\bullet\ra\P(\Delta)(\sigma)^\bullet$
is generated as a subcomplex by
\begin{equation}
\bigcup_{i=1}^{r_\sigma-1}\SdP(\Delta)(\sigma)^i
\end{equation}
for each $\sigma\in\Delta\setminus\{\Zero\}$ as we mensioned
before Lemma~\ref{Lem 1.5}.
Hence we get the surjection $\varphi$.
We get the lemma since
\begin{equation}
\Ker\varphi(\sigma)^{r_\sigma}\simeq
\Ker\varphi(\Delta,\t)(\sigma)^{r_\sigma}/
\Ker\psi_\Delta(\sigma)^{r_\sigma}
\end{equation}
and this is equal to
\begin{equation}
\bigoplus_{j=-r_\sigma}^{-1}
\psi_\Delta(\sigma)(\SdP(\Delta)(\sigma)_j^{r_\sigma})
= \bigoplus_{j=-r_\sigma}^{-1}\P(\Delta)(\sigma)_j^{r_\sigma}
\end{equation}
as a graded $\Q$-vector space.
This is equal to $\det(\sigma)\otimes(N(\sigma)A(\sigma))$.
\QED

We denote by $\bar A(\sigma)$ the $A(\sigma)$-module
$A(\sigma)/N(\sigma)A(\sigma)$ of lengh one.
By this lemma, we identify $\ic_\t(\Delta)^\bullet$ with the
quotient complex of $\P(\Delta)^\bullet$ by $\Ker\varphi$.
Namely we have
\begin{equation}
		\label{top}
\ic_\t(\Delta)(\sigma)^i =\left\{
\begin{array}{ll}
  \det(\sigma)\otimes\bar A(\sigma) &
  \hbox{ if }i = r_\sigma \\
  \{0\} & \hbox{ if }i\not= r_\sigma
\end{array}
\right.
\end{equation}
for every $\sigma\in\Delta$.


\begin{Lem}
		\label{Lem 1.10}
Let $\Delta$ be a finite fan.
Then $\H^p(\Gamma(\ic_\t(\Delta))^\bullet)_q =\{0\}$ for
$p, q\in\Z$ with $p > q + r$.
\end{Lem}

\Proof
Since
$\bar A(\sigma)_A =\bigwedge^\bullet(N_\Q/N(\sigma)_\Q)$
and $\dim_\Q N_\Q/N(\sigma)_\Q = r - r_\sigma$,
\begin{equation}
\Gamma(\ic_\t(\Delta))_q^p =
\bigoplus_{\sigma\in\Delta(p)}
\det(\sigma)\otimes(\bar A(\sigma)_A)_q =\{0\}
\end{equation}
for $q < -(r - p)$, i.e, for $p > q + r$.
Hence, for the component of degree $q$ of the complex
$\Gamma(\ic_\t(\Delta))^\bullet$, the cohomologies vanish for
$p > q + r$.
\QED


\section{The intersection complex of the middle perversity}
\setcounter{equation}{0}

Let $\Delta$ be a finite fan of $N_\R$.
The middle perversity $\m\cl\Delta\setminus\{\Zero\}\ra\Z$
is defined by $\m(\sigma) := 0$ for every $\sigma$.

We denote the intersection complex $\ic_\m(\Delta)^\bullet$
simply by $\ic(\Delta)^\bullet$.
Since $-\m =\m$, the dual $\D(\ic(\Delta))^\bullet$ is
quasi-isomorphic to $\ic(\Delta)^\bullet$ by
\cite[Cor.2.12]{Ishida3}.

A finite fan $\Delta$ of $N_\R$ is said to be a
{\em lifted complete fan}, if it is a lifting of a complete
fan of an $(r-1)$-dimensional space (cf.\cite[\S2]{Ishida3}).
In particular, $\Delta$ is a lifted complete fan if
the support $|\Delta|$ is equal to the boundary of an
$r$-dimensional cone.


\begin{Thm}
		\label{Thm 2.1}
Let $\Delta$ be a lifted complete fan of $N_\R$.
Then for any integers $p, q\in\Z$, the equality
\begin{equation}
		\label{zzliftedcomplete}
\dim_\Q\H^p(\Gamma(\ic(\Delta))^\bullet)_{q} =
\dim_\Q\H^{r-1-p}(\Gamma(\ic(\Delta))^\bullet)_{-r-q}
\end{equation}
holds.
\end{Thm}

\Proof
Since $\D(\ic(\Delta))^\bullet$ is quasi-isomorphic to
$\ic(\Delta)^\bullet$ by \cite[Cor.2.12]{Ishida3},
this is a consequence of \cite[Prop.2.8]{Ishida3}.
\QED


\begin{Cor}
		\label{Cor 2.2}
Let $\pi$ be an $r$-dimensional cone.
Then, for $\Delta := F(\pi)\setminus\{\pi\}$, the equality
(\ref{zzliftedcomplete}) folds for any integers $p, q\in\Z$.
\end{Cor}


Let $\Delta$ and $\Delta'$ be finite fans of $N_\R$.
Then $\Delta'$ is said to be a {\em subdivision} of $\Delta$
if $|\Delta| = |\Delta'|$ and, for every $\sigma\in\Delta'$
there exists $\rho\in\Delta$ with $\sigma\subset\rho$
(cf.\cite[Cor.1.16]{Oda1}).
If $\Delta'$ is a subdivision of $\Delta$, then there exists
a unique map $f\cl\Delta'\ra\Delta$ such that
$\sigma\cap\relint f(\sigma)\not=\emptyset$ for each
$\sigma\in\Delta'$.
Actually, $f(\sigma)$ is defined to be the minimal cone in
$\Delta$ which contains $\sigma$.
The map $f$ is also called a subdivision.

Let $f\cl\Delta'\ra\Delta$ be a subdivision.
For $\rho\in\Delta$, we denote by $f^{-1}(\rho)$ the
subset $f^{-1}(\{\rho\})$ of $\Delta'$.
Clearly, $f^{-1}(\rho)$ is a locally star closed subset
of $\Delta'$.

For $L\in\GEM(\Delta')$, we define the {\em direct image}
$f_*L\in\GEM(\Delta)$ by
\begin{equation}
(f_*L)(\rho) :=
\bigoplus_{\sigma\in f^{-1}(\rho)}L(\sigma)_{A(\rho)}
\end{equation}
for each $\rho\in\Delta$.

For a homomorphism $g\cl L\ra K$ in $\GEM(\Delta')$,
the homomorphism $f_*(g)\cl f_*L\ra f_*K$ is defined
as follows.

For $\rho,\mu\in\Delta$ with $\rho\prec\mu$ and
$\sigma\in f^{-1}(\rho),\;\tau\in f^{-1}(\mu)$,
the $(\sigma,\tau)$-component of
$f_*(g)(\rho/\mu)$ is the $A(\rho)$-homomorphism
$L(\sigma)_{A(\rho)}\ra K(\tau)_{A(\mu)}$
induced by the $A(\sigma)$-homomoprhism
$g(\sigma/\tau)\cl L(\sigma)\ra K(\tau)$
if $\sigma\prec\tau$, otherwise it is defined to be the
zero map.

It is easy to see that $f_*$ is a covariant additive
functor from $\GEM(\Delta')$ to $\GEM(\Delta)$.
We denote also by $f_*$ the induced functor from
$\CGEM(\Delta')$ to $\CGEM(\Delta)$.

A sequence of homomorphisms $0\ra L\ra K\ra J\ra 0$ in
$\GEM(\Delta)$ is said to be a short exact sequence if
$0\ra L(\sigma)\ra K(\sigma)\ra J(\sigma)\ra 0$ is exact
for every $\sigma\in\Delta$.


\begin{Lem}
\label{Lem 2.3}
Let $f\cl\Delta'\ra\Delta$ be a subdivision of a finite fan.

(1)
If $g\cl L^\bullet\ra K^\bullet$ is a quasi-isomoprhism in
$\CGEM(\Delta')$, then direct image
$f_*(g)\cl f_*L^\bullet\ra f_*K^\bullet$ is also a
quasi-isomorphism.

(2)
If $0\ra L\lraw{g}K\lraw{h}J\ra 0$ is a short exact
sequence in $\GEM(\Delta')$, then
$0\ra f_*L\lraw{f_*(g)}f_*K\lraw{f_*(h)}f_*J\ra 0$
is a short exact sequence in $\GEM(\Delta)$.
\end{Lem}

\Proof
(1)
Let $\Phi$ be a maximal subfan of $\Delta'$ such that
$f_*(g|\Phi)\cl f_*(L|\Phi)^\bullet\ra f_*(K|\Phi)^\bullet$
is a quasi-isomorphism.
Suppose $\Phi\not=\Delta'$ and let $\tau$ be a minimal
element of $\Delta'\setminus\Phi$.
Let $\Phi' :=\Phi\cup\{\tau\}$ and $\rho := f(\tau)$.
We get a commutative diagram
\begin{equation}
\begin{array}{ccccccccc}
0 & \ra & L(\tau)_{A(\rho)}^\bullet &
\lra & f_*(L|\Phi')(\rho)^\bullet &
\lra & f_*(L|\Phi)(\rho)^\bullet & \ra & 0 \\
 & & \hbox{  }\downarrow{\phi_1} & & \hbox{  }\downarrow{\phi_2}
 & & \hbox{  }\downarrow{\phi_3} & & \\
0 & \ra & K(\tau)_{A(\rho)}^\bullet &
\lra & f_*(K|\Phi')(\rho)^\bullet &
\lra & f_*(K|\Phi)(\rho)^\bullet & \ra & 0 \\
\end{array}
..
\end{equation}

Among the vertical homomorphisms, $\phi_1$ is a
quasi-isomorphism since
$g(\tau/\tau)\cl L(\tau)^\bullet\ra K(\tau)^\bullet$ is
quasi-isomorphic and $A(\rho)$ is a free $A(\tau)$-module,
while $\phi_3$ is quasi-isomorphic by the assumption.
Hence $\phi_2 = f_*(g|\Phi')(\rho)$ is also quasi-isomorhic.
Since $f_*(L|\Phi')(\mu) = f_*(L|\Phi)(\mu)$ for $\mu\not=\rho$,
$f_*(g|\Phi')$ is a quasi-isomorphism in $\CGEM(\Delta)$.
This contradicts the maximality of $\Phi$.
Hence $\Phi =\Delta'$.

(2)
Let $\Phi$ be a maximal subfan of $\Delta'$ such that
\begin{equation}
0\ra f_*(L|\Phi)\lraw{f_*(g)}f_*(K|\Phi)\lraw{f_*(h)}f_*(J|\Phi)\ra 0
\end{equation}
is a short exact sequence.
Suppose $\Phi\not=\Delta'$, and let $\tau$ be a minimal
element of $\Delta'\setminus\Phi$ and $\rho := f(\tau)$.
Let $\Phi' :=\Phi\cup\{\tau\}$.
It is sufficient to show that
\begin{equation}
0\ra f_*(L|\Phi')\lraw{f_*(g)}f_*(K|\Phi')\lraw{f_*(h)}f_*(J|\Phi')\ra 0
\end{equation}
is a short exact sequence, since it contradicts the maximality
of $\Phi$.
It is enough to check it for $\rho$.
Since
$0\ra L(\tau)_{A(\rho)}\lra K(\tau)_{A(\rho)}\lra J(\tau)_{A(\rho)}\ra 0$
is exact by the assumption, the exactness of
\begin{equation}
0\ra f_*(L|\Phi')(\rho)\lraw{f_*(g)}f_*(K|\Phi')(\rho)
\lraw{f_*(h)}f_*(J|\Phi')(\rho)\ra 0
\end{equation}
follows from the nine lemma of the homology algebra.
\QED

Although $\GEM(\Delta)$ and $\CGEM(\Delta)$ are not in
general abelian categories, we say that the functor $f_*$ is
exact in the sense that the property (2) of the above lemma holds.
The dualizing functor $\D$ defined in \cite[\S2]{Ishida3}
is an exact contravariant functor from $\CGEM(\Delta)$ to
itself in this sense.
Actually, for $L^\bullet\in\CGEM(\Delta)$ and $\rho\in\Delta$,
$\D(L)(\rho)^\bullet$ is defined by the combination of
exact functors $\i_\rho^*$ and $\d_\rho$ (cf.\cite[\S2]{Ishida3}).

Let $\Delta$ be a finite fan.
For each $\rho\in\Delta\setminus\{\Zero\}$, we take a
rational point $a(\rho)$ in the relative interior of $\rho$.
For each
$\alpha = (\rho_1,\cdots,\rho_k)\in\Sd(\Delta\setmin\{\Zero\})$,
let $c(\alpha)$ be the simplicial cone
$\R_0a(\rho_1) +\cdots +\R_0a(\rho_k)$ if $k > 0$, and
let $c(\alpha) :=\Zero$ if $k = 0$.
Then
\begin{equation}
\Sigma :=\{c(\alpha)\scl\alpha\in\Sd(\Delta\setmin\{\Zero\})\}
\end{equation}
is a simplicial subdivision of $\Delta$.

We call $\Sigma$ a {\em barycentric subdivision} of $\Delta$.
Clearly, it is not unique for $\Delta$ since it depends on
the choice of the set
$\{a(\rho)\scl\rho\in\Delta\setminus\{\Zero\}\}$.

Let $f\cl\Sigma\ra\Delta$ be a barycentric subdivision of
$\Delta$.
For $\alpha\in\Sd(\Delta\setmin\{\Zero\})$,
$f(c(\alpha)) =\max(\alpha)$ by the definition.
We define an unmixed homomorphism
$\lambda_\Delta\cl\SdP(\Delta)^\bullet\ra f_*\P(\Sigma)^\bullet$
as follows.

For each $\alpha\in\Sd(\Delta\setminus\{\Zero\})$, let
$c(\alpha)$ be the corresponding cone in $\Sigma$ as above.
For each $\rho\in\Delta$ and integer $i$, we have
\begin{equation}
\SdP(\Delta)(\rho)^i =
\bigoplus_{\beta\in\Sd_{i-1}(F(\Zero,\rho))}%
\Z z(\beta)y(\rho)\otimes A(\rho)
\end{equation}
and
\begin{equation}
\begin{array}{lll}
f_*\P(\Sigma)(\rho)^i & = &
\bigoplus_{\sigma\in f^{-1}(\rho)(i)}%
\det(\sigma)\otimes A(\rho) \\
 & = &\bigoplus_{\beta\in\Sd_{i-1}(F(\Zero,\rho))}%
\det(c(\beta)+c(\rho))\otimes A(\rho)\;,
\end{array}
\end{equation}
where $c(\rho)$ is the one-dimensional cone generated
by $a(\rho)$.

For each
$\beta = (\sigma_1,\cdots,\sigma_{i-1})\in\Sd_{i-1}(F(\Zero,\rho))$,
let $z'(\beta,\rho)$ be the generator of
$\det(c(\beta)+c(\rho))\simeq\Z$ which has the same sign with
$a(\sigma_1)\wedge\cdots\wedge a(\sigma_{i-1})\wedge a(\rho)$
in $\det(c(\beta)+c(\rho))\otimes\Q$.
We define $\lambda_\Delta(\rho)^i$ to be the isomoprhism given by
$z(\beta)y(\rho)\otimes 1\mapsto z'(\beta,\rho)\otimes 1$
for all $\beta\in\Sd_{i-1}(F(\Zero,\rho))$.

It is easy to check the compatibility with the coboundary
maps.
Since $\lambda_\Delta(\rho)^i$'s are isomorphic, we get the
following lemma.


\begin{Lem}
		\label{Lem 2.4}
Let $f\cl\Sigma\ra\Delta$ be a barycentric subdivision of
$\Delta$.
Then the above unmixed homomorphism
$\lambda_\Delta\cl\SdP(\Delta)^\bullet\ra f_*\P(\Sigma)^\bullet$
is an isomorphism.
\end{Lem}

Let $\psi_\Sigma\cl\SdP(\Sigma)^\bullet\ra\P(\Sigma)^\bullet$
be the natural unmixed quasi-isomorphism.
We define an unmixed homomorphism
\begin{equation}
\phi_{\Sigma/\Delta}\cl
f_*\SdP(\Sigma)^\bullet\lra\SdP(\Delta)^\bullet
\end{equation}
by $\phi_{\Sigma/\Delta} :=
\lambda_\Delta^{-1}\cdot f_*(\psi_\Sigma)$.


\begin{Lem}
		\label{Lem 2.5}
Let $f\cl\Sigma\ra\Delta$ be a barycentric subdivision and
let $\p$ and $\q$ be perversities of $\Delta$ and $\Sigma$,
respectively.
Then $\phi_{\Sigma/\Delta}(f_*k_\q(\Sigma))$ is contained
in $k_\p(\Delta)$ if $\q(\sigma)\leq\p(f(\sigma))$ for
every $\sigma\in\Sigma\setminus\{\Zero\}$.
In this case, $\phi_{\Sigma/\Delta}$ induces a natural
homomorphism
\begin{equation}
f_*\ic_\q(\Sigma)^\bullet\lra\ic_\p(\Delta)^\bullet\;.
\end{equation}
\end{Lem}

\Proof
By definition, $k_\q(\Sigma)$ is the subcomplex of
$\SdP(\Sigma)^\bullet$ generated by
\begin{equation}
S :=\bigcup_{\sigma\in\Sigma\setminus\{\Zero\}}\;%
\bigcup_{i+j\leq\q(\sigma)}\SdP(\Sigma)(\sigma)_j^i\;.
\end{equation}
Hence it is sufficient to show that
\begin{equation}
\phi_{\Sigma/\Delta}(f(\sigma))(\SdP(\Sigma)(\sigma)_j^i)
\subset k_\p(\Delta)(f(\sigma))
\end{equation}
for all $\sigma\in\Sigma\setminus\{\Zero\}$ and $i, j$
with $i+j\leq\q(\sigma)$.

If $i < r_\sigma$, then $\SdP(\Sigma)(\sigma)_j^i$ is
mapped to zero.
Hence the inclusion is obvious in this case.

We consider the case $i = r_\sigma$.
Let $\sigma\in\Sigma$ and let
$\alpha = (\sigma_1,\cdots,\sigma_{r_\sigma})$
be the corresponding element of
$\Sd_{r_\sigma}(\Delta\setminus\{\Zero\})$.
We set $\rho := f(\sigma) =\sigma_{r_\sigma}$.
Recall that $\SdP(\Sigma)(\sigma)^{r_\sigma}$ is the free
$A(\sigma)$-module with the basis
$\{z(\beta)y(\sigma)\scl\beta\in\Sd_{r_\sigma-1}(F(\Zero,\sigma))\}$.
Since $\phi_{\Sigma/\Delta}(\rho)(z(\beta)y(\sigma)\otimes 1)
=\pm z(\alpha)\otimes 1$
for every $\beta\in\Sd_{r_\sigma-1}(F(\Zero,\sigma))$,
\begin{equation}
\phi_{\Sigma/\Delta}(\rho)(\Z z(\beta)y(\sigma)\otimes A(\sigma)_j)
= \Z z(\alpha)\otimes A(\sigma)_j
\subset\Z z(\alpha)\otimes A(\rho)_j
\end{equation}
for every $\beta\in\Sd_{r_\sigma-1}(F(\Zero,\sigma))$
and $j\leq\q(\sigma) - r_\sigma$.
Hence the image $\phi_{\Sigma/\Delta}(\rho)(S(\sigma))$ is equal to
$\bigcup_{j\leq\q(\sigma) - r_\sigma}\Z z(\alpha)\otimes A(\sigma)_j$
and is contained in $k_\p(\Delta)(\rho)^{r_\sigma}$ if
$\q(\sigma)\leq\p(\rho)$.
\QED

Clearly, the condition of the above lemma is satisfied if
$\p$ and $\q$ are the middle perversities.
In particular we get a natural unmixed homomoprhism
\begin{equation}
		\label{zzichom}
\delta_{\Sigma/\Delta}\cl f_*\ic(\Sigma)^\bullet\lra\ic(\Delta)^\bullet
\end{equation}
for a barycentric subdivision $f\cl\Sigma\ra\Delta$.

Let $f\cl\Delta'\ra\Delta$ be a subdivision and
let $L^\bullet$ be an object of $\CGEM(\Delta')$.
We define an unmixed homomorphism
\begin{equation}
\kappa(f, L)\cl f_*\D(L)^\bullet\lra\D(f_*L)^\bullet
\end{equation}
as follows.

For $\rho\in\Delta$ and $i\in\Z$, we have
\begin{eqnarray}
\label{dual's image}
f_*\D(L)(\rho)^i & = &
\bigoplus_{\tau\in f^{-1}(\rho)}%
\D(L)(\tau)_{A(\rho)}^i \\
\label{dual's image2}
 & = &
\bigoplus_{\tau\in f^{-1}(\rho)}\;\bigoplus_{\sigma\in F(\tau)}%
\det(\tau)\otimes\d_\sigma(L(\sigma)^{r_\tau-i})_{A(\rho)}\;,
\end{eqnarray}
while
\begin{eqnarray}
\label{image's dual}
\D(f_*L)(\rho)^i & = &
\bigoplus_{\eta\in F(\rho)}%
\det(\rho)\otimes\d_\eta(f_*L(\eta)^{r_\rho-i})_{A(\rho)} \\
 & = &
\bigoplus_{\eta\in F(\rho)}\;\bigoplus_{\sigma\in f^{-1}(\eta)}%
\det(\rho)\otimes\d_\sigma(L(\sigma)^{r_\rho-i})_{A(\rho)} \\
 & = &
\bigoplus_{\sigma\in f^{-1}(F(\rho))}%
\det(\rho)\otimes\d_\sigma(L(\sigma)^{r_\rho-i})_{A(\rho)}\;.
\end{eqnarray}
For $\tau\in f^{-1}(\rho)$ and $\sigma\in F(\tau)$, the
restriction of
\begin{equation}
\kappa(f, L)(\rho)^i\cl f_*\D(L)(\rho)^i\lra\D(f_*L)(\rho)^i
\end{equation}
to the component
$\det(\tau)\otimes\d_\sigma(L(\sigma)^{r_\tau-i})_{A(\rho)}$
is defined to be the zero map if $r_\tau < r_\rho$.
If $r_\tau = r_\rho$, then $N(\tau) = N(\rho)$ and
$\det(\tau) =\det(\rho)$.
In this case, the component is defined to be the identity map to
$\det(\rho)\otimes\d_\sigma(L(\sigma)^{r_\rho-i})_{A(\rho)}$.

The commutativity of the diagram
\begin{equation}
\begin{array}{ccc}
\makebox[70pt]{}f_*\D(L)(\rho)^i &
\mathop{\lra}\limits^{d^i(\rho/\rho')} &
 f_*\D(L)(\rho')^{i+1}\makebox[50pt]{} \\
{\kappa(f, L)(\rho)^i}\downarrow &
 &\downarrow{\kappa(f, L)(\rho')^{i+1}} \\
\makebox[70pt]{}\D(f_*L)(\rho)^i &
\mathop{\lra}\limits^{d^i(\rho/\rho')} &
 \D(f_*L)(\rho')^{i+1}\makebox[50pt]{}
\end{array}
\end{equation}
is checked by the definitions.
The only one nontrivial case is the commutativity for the component
$\det(\tau)\otimes\d_\sigma(L(\sigma)^{r_\tau-i})_{A(\rho)}$
of $f_*\D(L)(\rho)^i$ with $r_\tau = r_\rho - 1$.
Since $\kappa(f, L)(\rho)^i$ is a zero map on this
component, we have to show that the composite
$\kappa(f, L)(\rho')^{i+1}\cdot d_{f_*\D(L)}^i(\rho/\rho')$
is zero on it.
Since $r_\tau = r_\rho - 1$, there exist exactly two cones
$\tau_1,\tau_2$ in $f^{-1}(\rho)$ with $\tau\prec\tau_1,\tau_2$
and $r_{\tau_1} = r_{\tau_2} = r_\rho$.
We have $q'_{\tau/\tau_1} + q'_{\tau/\tau_2} = 0$ under the
identification $\det(\tau_1) =\det(\tau_2) =\det(\rho)$.
Hence the composite is zero on this component.

Thus we know that $\kappa(f, L)$ is an unmixed homomorphism
in $\CGEM(\Delta)$.


\begin{Prop}
		\label{Prop 2.6}
Let $f\cl\Delta'\ra\Delta$ be a subdivision and
let $L^\bullet$ be an object of $\CGEM(\Delta')$.
Then the unmixed homomorphism
\begin{equation}
\kappa(f, L)\cl f_*\D(L)^\bullet\lra\D(f_*L)^\bullet
\end{equation}
is quasi-isomorphic.
\end{Prop}

\Proof
We prove the proposition by induction on the total dimension
of $L^\bullet$.
The assertion is trivially true if $\dim_\Q L = 0$.
We assume that $\dim_\Q L > 0$.

Let $\sigma$ be a maximal element of $\Delta'$ such that
$L(\sigma)^\bullet$ is nontrivial.
Take the maximal integer $p$ such that $L(\sigma)^p\not=\{0\}$
and the minimal integer $q$ such that $L(\sigma)_q^p\not=\{0\}$.

We define an object $E_{\sigma,p,q}^\bullet$ of $\CGEM(\Delta')$
as follows.

We define $E_{\sigma,p,q}(\tau)^\bullet :=\{0\}$ for
$\tau\in\Delta'$ with $\tau\not=\sigma$ and
$E_{\sigma,p,q}(\sigma)^i :=\{0\}$ for $i\not= p$.
The graded $A(\sigma)$-module $E_{\sigma,p,q}(\sigma)^p$
is defined to be $\bar A(\sigma)(-q)$, i.e.,
\begin{equation}
E_{\sigma,p,q}(\sigma)_j^p :=\left\{
\begin{array}{ll}
\Q    & \hbox{ if } j = q \\
\{0\} & \hbox{ if } j\not= q
\end{array}
\right.\;.
\end{equation}

By taking a nontrivial $A(\sigma)$-homomoprhism
$g_0\cl E_{\sigma,p,q}(\sigma)^p\ra L(\sigma)^p$, we get an unmixed
homomorphism $g\cl E_{\sigma,p,q}^\bullet\ra L^\bullet$
such that $g(\sigma)^p = g_0$.
Let $K^\bullet$ be the cokernel of $g$.
Then we get a commutative diagram
\begin{equation}
\begin{array}{ccccccccc}
0 & \ra & f_*\D(K)^\bullet &
\lra & f_*\D(L)^\bullet &
\lra & f_*\D(E_{\sigma,p,q})^\bullet & \ra & 0 \\
 & & \hbox{  }\downarrow{\kappa(f, K)}
 & & \hbox{  }\downarrow{\kappa(f, L)}
 & & \hbox{  }\downarrow{\kappa(f, E_{\sigma,p,q})} & & \\
0 & \ra & \D(f_*K)^\bullet &
\lra & \D(f_*L)^\bullet &
\lra & \D(f_*E_{\sigma,p,q})^\bullet & \ra & 0 \\
\end{array}
\;.
\end{equation}

Since the functors $f_*$ and $\D$ are exact, the two horizontal
lines in the diagram are short exact sequences.
Since $\dim_\Q K =\dim_\Q L - 1$, $\kappa(f, K)$ is a
quasi-isomorphism by the induction assumption.
Hence it is sufficient to show that $\kappa(f, E_{\sigma,p,q})$
is quasi-isomorphic.

Let $\rho$ be an element of $\Delta$.
If $f(\sigma)$ is not in $F(\rho)$, then both
$f_*\D(E_{\sigma,p,q})(\rho)^\bullet$ and
$\D(f_*E_{\sigma,p,q})(\rho)^\bullet$ are zero.
We assume $\eta := f(\sigma)\in F(\rho)$.
Set
\begin{equation}
\Phi :=\{\tau\in\Delta'\scl\sigma\prec\tau\in f^{-1}(\rho)\}\;.
\end{equation}
By (\ref{dual's image2}), we have
\begin{equation}
f_*\D(E_{\sigma,p,q})(\rho)^i =
\bigoplus_{\tau\in\Phi(p+i)}\det(\tau)\otimes\bar A(\sigma)_{A(\rho)}\;.
\end{equation}
Hence we know
\begin{equation}
\H^i(f_*\D(E_{\sigma,p,q})(\rho)^\bullet)\simeq
\H^{p+i}(E(\Phi,\Z)^\bullet)\otimes\bar A(\sigma)_{A(\rho)}
\end{equation}
for $i\in\Z$.
On the other hand, $\D(f_*E_{\sigma,p,q})(\rho)^i$ is zero
if $i\not= r_\rho - p$, while
\begin{equation}
\D(f_*E_{\sigma,p,q})(\rho)^{r_\rho - p} =
\det(\rho)\otimes\bar A(\sigma)_{A(\rho)}
\end{equation}
by (\ref{image's dual}).

The cohomologies $\H^i(E(\Phi,\Z)^\bullet)$ are zero for
$i\not= r_\rho$ and
$\H^{r_\rho}(E(\Phi,\Z)^\bullet) =\det(\rho)$ by
Lemma~\ref{Lem 1.2},(2).
We know $\H^i(f_*\D(E_{\sigma,p,q})(\rho)^\bullet)$ is
zero for $i\not= r_\rho - p$ and isomorphic to
$\det(\rho)\otimes\bar A(\sigma)_{A(\rho)}$ for $i = r_\rho - p$.
Since $\kappa(f, E_{\sigma,p,q})(\rho)^{r_\rho-p}$ is surjective,
$\kappa(f, E_{\sigma,p,q})(\rho)$ is quasi-isomorphic.
\QED

The intersection complex has the following irreducibility.


\begin{Lem}
		\label{Lem 2.7}
Let $\Delta$ be a finite fan, and $\p$ a perversity
of it.
Let $ L_1^\bullet$ and $L_2^\bullet$ be objects of
$\CGEM(\Delta)$ which are quasi-isomorphic to
$\ic_\p(\Delta)$.
If $h\cl L_1^\bullet\ra L_2^\bullet$ is a homomoprhism
such that
$h(\Zero/\Zero)\cl L_1(\Zero)^\bullet\ra L_2(\Zero)^\bullet$
is a quasi-isomorphism, then $h$ is a quasi-isomorphism.
\end{Lem}

\Proof
Suppose that $h$ is not quasi-isomorphic.
Let $\rho$ be a minimal element of $\Delta$ such that
$h(\rho/\rho)\cl L_1(\rho)^\bullet\ra L_2(\rho)^\bullet$
is not quasi-isomorphic.
By the assumption, $\rho$ is not equal to $\Zero$.
By the minimality of $\rho$, the homomorphism
\begin{equation}
\i_\rho^\circ(h)\cl
\i_\rho^\circ(L_1)^\bullet\ra\i_\rho^\circ(L_2)^\bullet
\end{equation}
in $\CGM(A(\rho))$ is a quasi-isomorphism.

For each of integer $i$, we have a
commutative diagram
\begin{equation}
\begin{array}{ccc}
\makebox[10pt]{}\H^i(\i_\rho^\circ L_1^\bullet) &
\mathop{\lra}\limits^{u^i} &
 \H^i(\i_\rho^\circ L_2^\bullet)\makebox[10pt]{} \\
{\phi_1}\downarrow &
 &\downarrow{\phi_2} \\
\makebox[10pt]{}\H^i(L_1(\rho)[1]^\bullet) &
\mathop{\lra}\limits^{v^i} &
 \H^i(L_2(\rho)[1]^\bullet)\makebox[10pt]{}
\end{array}
,
\end{equation}
where $u^i$ and $v^i$ are the $A(\rho)$-homomorphism induced
by $\i_\rho^\circ(h)$ and $h(\rho/\rho)$, respectively.
Since $\i_\rho^\circ(h)$ is quasi-isomorphic, $u^i$ is an
isomorphism.
Since $L_2^\bullet$ is quasi-isomorphic to $\ic_\p(\Delta)^\bullet$,
$\phi_2$ is surjective by the construction of intersection
complexes \cite[Thm.2.9]{Ishida3}.
Hence $v^i$ is also surjective.
Since $\H^i(L_1(\rho)[1]^\bullet)$ and
$\H^i(L_2(\rho)[1]^\bullet)$ are finite dimensional
$\Q$-vector spaces of same dimension, $v^i$ is an
isomorphism for each $i$, i.e., $h(\rho/\rho)$ is a
quasi-isomorphism.
This contradicts the assumption.
\QED

Since $\ic_\p(\Delta)(\Zero)^0 =\Q$ and
$\ic_\p(\Delta)(\Zero)^i =\{0\}$ for $i\not= 0$,
$h(\Zero/\Zero)$ in the above lemma is quasi-isomorphic
if and only if the induced homomorphism
$\H^0(L_1(\Zero)^\bullet)\ra\H^0(L_2(\Zero)^\bullet)$
is a nonzero map.

Let $f\cl\Delta'\ra\Delta$ be a subdivision and let
$L^\bullet$ be in $\CGEM(\Delta')$.
Then, by the definition of the functor $f_*$,
the complex $\Gamma(f_*L)^\bullet$ is canonically
isomorphic to $\Gamma(L)^\bullet$.
In particular, if $f\cl\Sigma\ra\Delta$ is a barycentric
subdivision, then $\delta_{\Sigma/\Delta}$ in
(\ref{zzichom}) induces a homomorphism
\begin{equation}
\Gamma(\delta_{\Sigma/\Delta})\cl\Gamma(\ic(\Sigma))^\bullet\lra
\Gamma(\ic(\Delta))^\bullet
\end{equation}
in $\CGM(A)$.


\begin{Thm}[Decomposition theorem]
		\label{Thm 2.8}
Let $\Delta$ be a finite fan and $\Sigma$ a barycentric
subdivision of $\Delta$.
Then, for each integer $p$, the homomorphism of $A$-modules
\begin{equation}
		\label{zzhomThm 2.7}
\H^p(\Gamma(\ic(\Sigma))^\bullet)\lra
\H^p(\Gamma(\ic(\Delta))^\bullet)
\end{equation}
induced by $\Gamma(\delta_{\Sigma/\Delta})$ is a split surjection, i.e.,
is surjective and the kernel is a direct summand as an
$A$-module.
In particular,
\begin{equation}
\dim_\Q\H^p(\Gamma(\ic(\Delta))^\bullet)_{q}\leq
\dim_\Q\H^p(\Gamma(\ic(\Sigma))^\bullet)_{q}
\end{equation}
for any integers $p$, $q$.
\end{Thm}

\Proof
By applying the contravariant functor $\D$ to the
homomorphism (\ref{zzichom}), we get a homomorphism
\begin{equation}
		\label{zzrevichom}
\D(\delta_{\Sigma/\Delta})\cl\D(\ic(\Delta))^\bullet
\lra\D(f_*\ic(\Sigma))^\bullet\;.
\end{equation}
Since $\delta_{\Sigma/\Delta}(\Zero)$ is an isomorphism, so is
$\D(\delta_{\Sigma/\Delta})(\Zero)$.

Since $\D(\ic(\Sigma))^\bullet$ is quasi-isomorphic to
$\ic(\Sigma)^\bullet$ by \cite[Cor.2.12]{Ishida3},
$f_*\D(\ic(\Sigma))^\bullet$ and $f_*\ic(\Sigma)^\bullet$
are also quasi-isomorphic by Lemma~\ref{Lem 2.3},(1).
By Proposition~\ref{Prop 2.6}, we get a quasi-isomorphism
\begin{equation}
\kappa(f,\ic(\Sigma))\cl
f_*\D(\ic(\Sigma))^\bullet\lra\D(f_*\ic(\Sigma))^\bullet\;.
\end{equation}
Hence $\D(f_*\ic(\Sigma))^\bullet$ is quasi-isomorphic to
$f_*\ic(\Sigma)^\bullet$.

On the other hand, $\D(\ic(\Delta))^\bullet$ is quasi-isomorphic
to $\ic(\Delta)^\bullet$ by \cite[Cor.2.12]{Ishida3}.

By applying \cite[Lem 2.16]{Ishida3} for the homomorphism
$\D(\delta_{\Sigma/\Delta})$, we get $L^\bullet$ in $\CGEM(\Delta)$,
a quasi-isomorphism $g\cl L^\bullet\ra\ic(\Delta)^\bullet$
and a homomorphism
$h\cl L^\bullet\ra f_*\ic(\Sigma)^\bullet$ such that the
homomorphisms of cohomologies induced by $h$ is compatible
with those induced by $\D(\delta_{\Sigma/\Delta})$.

Note that all homomorphisms in $\CGEM(\Delta)$ appeared
here are quasi-isomorphic on $\Zero\in\Delta$.
Since $L^\bullet$ and $\ic(\Delta)^\bullet$ are quasi-isomorphic,
the composite $\delta_{\Sigma/\Delta}\cdot h\cl
L^\bullet\ra\ic(\Delta)^\bullet$
is a quasi-isomorphism by Lemma~\ref{Lem 2.7}.

Hence, for each integer $i$, the homomorphism (\ref{zzhomThm 2.7})
is surjective and $\H^p(\Gamma(\ic(\Sigma))^\bullet)$ is the
direct sum of the kernel and the image of the homomorphism
$\H^p(L^\bullet)\ra\H^p(\Gamma(\ic(\Sigma))^\bullet)$ induced
by $h$.
\QED


\section{Intersection complexes for simplicial fans}
\setcounter{equation}{0}


\begin{Lem}
		\label{Lem 3.1}
Let $\pi$ be a simplicial cone.
Then
\begin{equation}
\H^i(\ic_\t(F(\pi))(\pi)^\bullet)_j\simeq\left\{
\begin{array}{ll}
\Q    & \hbox{ if } (i, j) = (r_\pi, 0) \\
\{0\} & \hbox{ if } (i, j)\not= (r_\pi, 0)
\end{array}
\right.\;.
\end{equation}
On the other hand,
\begin{equation}
\H^i(\i_\pi^*\ic_\t(F(\pi))^\bullet)_j\simeq\left\{
\begin{array}{ll}
\Q    & \hbox{ if } (i, j) = (0, -r_\pi) \\
\{0\} & \hbox{ if } (i, j)\not= (0, -r_\pi)
\end{array}
\right.\;.
\end{equation}
\end{Lem}

\Proof
The first assertion follows from the description (\ref{top}).

Let $s := r_\pi$ and $\{x_1,\cdots, x_s\}\subset N_\Q$ be a
minimal generator of the simplicial cone $\pi$.
For each element $\sigma\in F(\pi)$, there exists a unique
subset $\{i_1,\cdots, i_p\}$ of $\{1,\cdots, s\}$ with
$i_1 <\cdots < i_p$ such that $\{x_{i_1},\cdots, x_{i_s}\}$
generates $\sigma$.
We denote $x(\sigma) := x_{i_1}\wedge\cdots\wedge x_{i_s}$.

For $\rho\in F(\pi)$, we have a description
\begin{equation}
\ic_\t(F(\pi))(\rho)^i_{A(\pi)} =
\det(\rho)\otimes{\bigwedge}\!^\bullet(N(\pi)_\Q/N(\rho)_\Q) =
\bigoplus_{\sigma\in F(\rho')}\det(\rho)\otimes\Q x(\sigma)\;,
\end{equation}
where $\rho'$ is the complementary cone of $\rho$ in $\pi$,
i.e., the unique cone in $F(\pi)$ with $\rho\cap\rho' =\Zero$
and $\rho +\rho' =\pi$.
Hence $\ic_\t(F(\pi))(\rho)_{A(\pi)}^i$ has the
component $\det(\rho)\otimes\Q x(\sigma)$ if and only if
$\rho\in F(\sigma')$ for the complementary cone $\sigma'$
of $\sigma$.
Note that $x(\sigma)$ is a homogeneous element of degree
$-r_\sigma$ and $\sigma'$ is of dimension $r_\pi - r_\sigma$.

By this observation, we have
\begin{equation}
(\i_\pi^*\ic_\t(F(\pi))^\bullet)_j\simeq
\bigoplus_{\sigma\in F(\pi)(-j)}E(F(\sigma'),\Z)^\bullet\otimes\Q x(\sigma)
\end{equation}
for each integer $j$.
By Lemma~\ref{Lem 1.2},(1), this complex of $\Q$-vector spaces
has trivial cohomologies if $r_\pi + j > 0$ since
$\dim\sigma' = r_\pi + j$ for $\sigma\in F(\pi)(-j)$.

For $j = - r_\pi$, $F(\pi)(-j) =\{\pi\}$ and $\pi' =\Zero$.
Hence $\H^0(\i_\pi^*\ic_\t(F(\pi))^\bullet)_{-r_\pi}\simeq\Q$
and $\H^i(\i_\pi^*\ic_\t(F(\pi))^\bullet)_{-r_\pi} =\{0\}$
for $i\not= 0$.

We get the lemma since the complexes are nontrivial only
for $-r_\pi\leq j\leq 0$.
\QED


\begin{Thm}
		\label{Thm 3.2}
Let $\Delta$ be a simplicial finite fan.
Then, for any perversity $\p$ with $\b\leq\p\leq\t$,
$\ic_\p(\Delta)^\bullet$ is quasi-isomorphic to
$\ic_\t(\Delta)^\bullet$.
\end{Thm}

\Proof
By \cite[Thm.2.9]{Ishida3}, it sufficient to show that
$\H^i(\ic_\p(\Delta)(\rho)^\bullet)_j =\{0\}$
if $i +j\leq\p(\rho)$ and
$\H^i(\i_\rho^*(\ic_\p(\Delta))^\bullet)_j =\{0\}$
if $i +j\geq\p(\rho)$ for $\rho\in\Delta\setminus\{\Zero\}$.
Since $\ic_\p(F(\rho))^\bullet$ is the restiction of
$\ic_\p(\Delta)^\bullet$ to $F(\rho)$ and
$-r_\rho + 1\leq\p(\rho)\leq r_\rho - 1$ by the assumption, this
is a consequence of Lemma~\ref{Lem 3.1}.

The following theorem is equivalent to \cite[Prop.4.1]{Oda2}.


\begin{Thm}
		\label{Thm 3.3}
Let $\Delta$ be a simplicial complete fan.
Then $\H^p(\Gamma(\ic_\t(\Delta))^\bullet)_q =\{0\}$ for
any integers $p, q\in\Z$ with $p\not= q + r$.
\end{Thm}

\Proof
By Lemma~\ref{Lem 1.10}, we have
$\H^p(\Gamma(\ic_\t(\Delta))^\bullet)_q =\{0\}$ for
$p > q + r$.
By the duality for complete fans \cite[Prop.2.5]{Ishida3},
we have $\H^p(\Gamma(\D(\ic_\t(\Delta)))^\bullet)_q =\{0\}$
for $p, q\in\Z$ with $(r-p) > (-r-q) + r$, i.e.,
with $p < q + r$.
$\D(\ic_\t(\Delta))^\bullet$ is quasi-isomorphic to
$\ic_\b(\Delta)^\bullet$ by \cite[Cor.2.12]{Ishida3}.
We get the theorem, since $\ic_\b(\Delta)^\bullet$ is
quasi-isomorphic to $\ic_\t(\Delta)^\bullet$ by
Theorem~\ref{Thm 3.2}.
\QED

Let $f\cl\Delta'\ra\Delta$ be an arbitrary subdivision of
a finite fan.
We define an unmixed homomorphism
\begin{equation}
\label{image of top}
\bar\delta_{\Delta'/\Delta}\cl
f_*\ic_\t(\Delta')^\bullet\ra\ic_\t(\Delta)^\bullet
\end{equation}
in $\CGEM(\Delta)$ as follows.

Note that
\begin{equation}
f_*\ic_\t(\Delta')(\rho)^i :=
\bigoplus_{\sigma\in f^{-1}(\rho)(i)}
\det(\sigma)\otimes(\bar A(\sigma))_{A(\rho)}
\end{equation}
for each integer $i$, while
\begin{equation}
\ic_\t(\Delta)(\rho)^i = \left\{
\begin{array}{ll}
  \det(\rho)\otimes\bar A(\rho) &
  \hbox{ if }i = r_\rho\\
  \{0\} & \hbox{ if }i\not= r_\rho
\end{array}
\right.\;.
\end{equation}

For each $\sigma\in f^{-1}(\rho)$, the restriction to the
component $\det(\sigma)\otimes(\bar A(\sigma))_{A(\rho)}$
of the homomorphism
\begin{equation}
\bar\delta_{\Delta'/\Delta}(\rho)^i\cl
f_*\ic_\t(\Delta')(\rho)^i\ra\ic_\t(\Delta)(\rho)^i
\end{equation}
is defined to be zero if $i = r_\sigma < r_\rho$.
If $r_\sigma = r_\rho$, then $\det(\sigma) =\det(\rho)$
and $A(\sigma) = A(\rho)$, and the component is defined to
be the identity map to $\det(\rho)\otimes\bar A(\rho)$.

Let $N'$ be a free $\Z$-module of rank $r'$ and let
$f_0\cl N\ra N'$ be a surjective homomorphism.
Fans $\Delta$ of $N_\R$ and $\Phi$ of $N'_\R$ are said to
be compatible with $f_0$, if $f_0(\sigma)$ is contained in
a cone of $\Phi$ for every $\sigma\in\Delta$, where we
denote also by $f_0$ the linear map $N_\R\ra N'_\R$.

We say a map $f\cl\Delta\ra\Phi$ a {\em morphism of fans}
if such $f_0$ is given and $f(\sigma)$ is the minimal cone
of $\Phi$ which contains $f_0(\sigma)$ for every
$\sigma\in\Delta$.
The direct image functor
\begin{equation}
f_*\cl\GEM(\Delta)\lra\GEM(\Phi)
\end{equation}
is defined as follows.

Let $L$ be in $\GEM(\Delta)$.
For each $\rho\in\Phi$, we set $f^{-1}(\rho) := f^{-1}(\{\rho\})$
similarly as in the case of subdivisions.
If $\sigma\in f^{-1}(\rho)$, then $f_0$ induces a homomorphism
$N(\sigma)\ra N'(\rho)$ and a ring homomorphism
$A(\sigma)\ra A'(\rho)$.
We set
\begin{equation}
(f_*L)(\rho)  :=\bigoplus_{\sigma\in f^{-1}(\rho)}
L(\sigma)_{A'(\rho)}\;,
\end{equation}
where $L(\sigma)_{A'(\rho)} := L(\sigma)\otimes_{A(\sigma)}A'(\rho)$.
Let $g\cl L\ra K$ be a homomorphism in $\GEM(\Delta)$.
For $\rho,\mu\in\Phi$ with $\rho\prec\mu$ and for
$\sigma\in f^{-1}(\rho)$ and $\tau\in f^{-1}(\mu)$, the
$(\sigma,\tau)$-component of the homomorphism
$f_*(g)(\rho/\mu)\cl(f_*L)(\rho)\ra(f_*K)(\mu)$
is the map induced by $g(\sigma/\tau)$ if $\sigma\prec\tau$
and is the zero map otherwise.

Let $\Delta$ and $\Phi$ be complete fans of $N_\R$ and
$N'_\R$, respectively, and let $f\cl\Delta\ra\Phi$ be a
morphism of fans associated with $f_0\cl N\ra N'$.
We take a generator $w$ of $\bigwedge^{r-r'}\ker f_0\simeq\Z$.
We define an unmixed homomorphism
\begin{equation}
\label{from direct image}
\bar\delta_{\Delta/\Phi}\cl
f_*\ic_\t(\Delta)^\bullet\lra\ic_\t(\Phi)[-r+r']^\bullet
\end{equation}
as follows.

Let $\rho$ be in $\Phi$.
For each $\sigma\in f^{-1}(\rho)$, the
$\sigma$-component of the homomorphism
\begin{equation}
\bar\delta_{\Delta/\Phi}(\rho)^i\cl f_*\ic_\t(\Delta)(\rho)^i
\lra\ic_\t(\Phi)[-r+r'](\rho)^i
\end{equation}
is defined to be zero if $i\not= r_\rho + r -r'$ since
then
\begin{equation}
\ic_\t(\Phi)[-r+r'](\rho)^i =
\ic_\t(\Phi)(\rho)^{i-r+r'} =\{0\}\;.
\end{equation}
If $i = r_\rho + r -r'$, then
\begin{equation}
f_*\ic_\t(\Delta)(\rho)^i =
\bigoplus_{\sigma\in f^{-1}(\rho)(r_\rho + r -r')}
\det(\sigma)\otimes\bar A(\sigma)_{A'(\rho)}
\end{equation}
and
\begin{equation}
\ic_\t(\Phi)[-r+r'](\rho)^i =
\det(\rho)\otimes\bar A'(\rho)\;.
\end{equation}
If $\sigma\in f^{-1}(\rho)(r_\rho + r -r')$, then
$r_\sigma - r_\rho = r - r'$ and $\Ker f_0\subset N(\sigma)$.
Hence $w\in\bigwedge^{r-r'}N(\sigma)$.
The $\sigma$-component of the homomorphism is defined to
be the tensor product of the isomorphism
$q_{\sigma/\rho}^w\cl\det(\sigma)\ra\det(\rho)$ which sends
$w\wedge a$ to $(\bigwedge^{r_\rho}f_0)(a)$ for
$a\in\bigwedge^{r_\rho}N(\sigma)$ and the natural isomorphism
$\bar A(\sigma)_{A'(\rho)}\simeq\bar A'(\rho)$, where
$\bigwedge^{r_\rho}f_0\cl\bigwedge^{r_\rho}N(\sigma)\ra\det(\rho)$
is the natural map induced by $f_0$.

Let $\eta$ be an element of a finite fan $\Delta$ of $N_\R$
and let $N'$ be the quotient free $\Z$-module $N/N(\eta)$ of
rank $r' := r - r_\eta$.
We set $\Delta(\eta{\prec}) :=
\{\sigma\in\Delta\scl\eta\prec\sigma\}$ as before.
For each $\sigma\in\Delta(\eta{\prec})$, let $\sigma[\eta]$ be
the image of $\sigma$ in the quotient space
$N'_\R = N_\R/N(\eta)_\R$.
Then
\begin{equation}
\Delta[\eta] :=\{\sigma[\eta]\scl\sigma\in\Delta(\eta{\prec})\}
\end{equation}
is a finite fan of $N'_\R$.
For $\tau\in\Delta[\eta]$, the notations $N'(\tau)\subset N'$ and
$A'(\tau)\subset A' :=\bigwedge^{r - r_\eta}N'_\Q$ are defined
similarly for $N'$ as we defined for $N$.

For $\sigma\in\Delta(\eta{\prec})$, $A'(\sigma[\eta])$
is a quotient ring of $A(\sigma)$ by the two-sided ideal
generated by $N(\eta)$.
Hence the category $\GM(A'(\sigma[\eta]))$ of fintely generated
graded $A'(\sigma[\eta])$-modules is a full subcategory of
$\GM(A(\sigma))$.

We define the categories $\GEM(\Delta[\eta])$ and
$\CGEM(\Delta[\eta])$, similarly.
We denote by $\epsilon_\eta$ the natural functor
\begin{equation}
\epsilon_\eta\cl\GEM(\Delta[\eta])\lra\GEM(\Delta)
\end{equation}
defined by $\epsilon_\eta(L)(\sigma) := L(\sigma[\eta])$
if $\sigma\in\Delta(\eta{\prec})$ and
$\epsilon_\eta(L)(\sigma) :=\{0\}$ otherwise.
For $f\cl L\ra K$ in $\GEM(\Delta[\eta])$,
$\epsilon_\eta(f)\cl
\epsilon_\eta(L)\ra\epsilon_\eta(K)$ is defined by
$\epsilon_\eta(f)(\sigma/\tau) := f(\sigma[\eta]/\tau[\eta])$
if $\sigma\in\Delta(\eta{\prec})$ and
$\epsilon_\eta(f)(\sigma/\tau) := 0$ otherwise.

We take a generator $w\in\det(\eta) =\Ker f_0\simeq\Z$.
We define an isomorphism
\begin{equation}
\label{quotient iso}
h(\Delta,\eta, w)\cl
\ic_\t(\Delta(\eta{\prec}))^\bullet\simeq
\epsilon_\eta(\ic_\t(\Delta[\eta]))[-r_\eta]^\bullet
\end{equation}
as follows.

For each $\sigma\in\Delta(\eta{\prec})$, we have
\begin{equation}
\ic_\t(\Delta(\eta{\prec}))(\sigma)^{r_\sigma} =
\det(\sigma)\otimes\bar A(\sigma)\;,
\end{equation}
and
\begin{equation}
\epsilon_\eta(\ic_\t(\Delta[\eta]))[-r_\eta](\sigma)^{r_\sigma}
= \det(\sigma[\eta])\otimes\bar A'(\sigma[\eta])
\end{equation}
by (\ref{top}), while
\begin{equation}
\ic_\t(\Delta(\eta{\prec}))(\sigma)^i =
\epsilon_\eta(\ic_\t(\Delta[\eta]))[-r_\eta](\sigma)^i =\{0\}
\end{equation}
for $i\not = r_\sigma$.
Since $\bar A(\sigma) = A(\sigma)/N(\sigma)A(\sigma)$ and
$\bar A'(\sigma[\eta]) =
A'(\sigma[\eta])/N'(\sigma[\eta])A'(\sigma[\eta])$
and since the kernel of the surjection
$A(\sigma)\ra A'(\sigma[\eta])$ is
$N(\eta)A(\sigma)\subset N(\sigma)A(\sigma)$,
$\bar A'(\sigma[\eta])$ is naturally isomorphic to
$\bar A(\sigma)$ as an $A(\sigma)$-module.
we define an isomorphism
\begin{equation}
h(\Delta,\eta, w)(\sigma)^{r_\sigma}\cl
\ic_\t(\Delta(\eta{\prec}))(\sigma)^{r_\sigma}\lra
\epsilon_\eta(\ic_\t(\Delta[\eta]))[-1](\sigma)^{r_\sigma}
\end{equation}
to be the tensor product of $q_{\sigma/\sigma[\eta]}^w$ and
this isomorphism.
The commutativity with the coboundary maps is checked
easily.
When $r_\eta = 1$, we take the primitive element of
$\eta\cap N$ as $w$ and denote the ismorphism simply by
$h(\Delta,\eta)$.


\begin{Lem}
\label{Lem 3.4}
Let $\eta$ be a cone of a complete fan $\Delta$.
We define $N'$, $f_0\cl N_\R\ra N'_\R$ and $\Delta[\eta]$
as above.
If $f_0$ induces a morphism $f\cl\Delta\ra\Delta[\eta]$ of
fans, then the homomorhism
\begin{equation}
\H^i(\Gamma(\ic_\t(\Delta(\eta{\prec})))^\bullet)\lra
\H^i(\Gamma(\ic_\t(\Delta))^\bullet)
\end{equation}
induced by the inclusion map
$\ic_\t(\Delta(\eta{\prec}))^\bullet\ra\ic_\t(\Delta)^\bullet$
is injective for every $i\in\Z$.
\end{Lem}

\Proof
We take a generator $w$ of $\det(\eta)$.
By definitions, the composite of homomorphisms
\begin{equation}
\label{composite 3.4}
\Gamma(\ic_\t(\Delta(\eta{\prec})))^\bullet
\mathop{\lra}\limits^{\phi}
\Gamma(\ic_\t(\Delta))^\bullet\ra
\Gamma(f_*\ic_\t(\Delta))^\bullet\ra
\Gamma(\ic_\t(\Delta[\eta])[-r_\eta])^\bullet
\end{equation}
induced by the homomorphisms
$\ic_\t(\Delta(\eta{\prec}))^\bullet\ra
\ic_\t(\Delta)^\bullet$
and
$\bar\delta_{\Delta/\Delta[\eta]}\cl f_*\ic_\t(\Delta)^\bullet\ra
\ic_\t(\Delta[\eta])[-r_\eta]^\bullet$
is equal to the homomoprhism
\begin{equation}
\Gamma(\ic_\t(\Delta(\eta{\prec})))^\bullet\lra
\Gamma(\ic_\t(\Delta[\eta])[-r_\eta])^\bullet
\end{equation}
induced by $h(\Delta,\eta, w)$.
Since $h(\Delta,\eta, w)$ is an isomorphism, the
homomorphisms of cohomologies induced by $\phi$ in
(\ref{composite 3.4}) are injective.
\QED


\begin{Thm}
\label{Thm 3.5}
Let $\Delta$ be a simplicial complete fan and let
$\eta$ be an element of $\Delta$.
Then the homomorphism
\begin{equation}
\H^i(\Gamma(\ic_\t(\Delta(\eta{\prec})))^\bullet)\lra
\H^i(\Gamma(\ic_\t(\Delta))^\bullet)
\end{equation}
induced by the inclusion map
$\ic_\t(\Delta(\eta{\prec}))^\bullet\ra\ic_\t(\Delta)^\bullet$
is injective for every $i\in\Z$.
\end{Thm}

\Proof
Let $g\cl L^\bullet\ra K^\bullet$ be a quasi-isomorphism in
$\CGEM(\Delta)$.
Then we get a commutative diagram
\begin{equation}
\begin{array}{ccc}
\makebox[20pt]{}\H^i(\Gamma(L|\Delta(\eta{\prec}))^\bullet) &
\mathop{\lra}\limits^{\phi_1} &
 \makebox[10pt]{}\H^i(\Gamma(L)^\bullet) \\
{g_1}\downarrow &
 &\downarrow{g_2} \\
\makebox[20pt]{}\H^i(\Gamma(K|\Delta(\eta{\prec}))^\bullet) &
\mathop{\lra}\limits^{\phi_2} &
\makebox[10pt]{}\H^i(\Gamma(K)^\bullet)
\end{array}
,
\end{equation}
for every $i\in\Z$.
Since $g$ is a quasi-isomorphism, $g_1$ and $g_2$ are isomorphisms.
Hence $\phi_1$ is injective if and only if $\phi_2$ is
injective.
This implies that, for the proof of the theorem, it is sufficient
to show the injectivity of the homomorphisms
\begin{equation}
\H^i(\Gamma(L|\Delta(\eta{\prec}))^\bullet)\lra
\H^i(\Gamma(L)^\bullet)
\end{equation}
for $i\in\Z$ for an $L^\bullet$ which is quasi-isomorphic to
$\ic_\t(\Delta)^\bullet$.

Let $N' := N/N(\eta)$ and $f_0\cl N_\R\ra N'_\R$ be the
natural surjection.
We take a simplicial subdivision $u\cl\Delta'\ra\Delta$ such
that $\eta\in\Delta'$,
$\Delta'(\eta{\prec}) =\Delta(\eta{\prec})$ and
$f_0$ induces a morphism $f\cl\Delta'\ra\Delta'[\eta]$.

Then $u_*\ic_\t(\Delta'(\eta{\prec}))^\bullet$ is
equal to $\ic_\t(\Delta(\eta{\prec}))^\bullet$.
Since
$\Gamma(u_*\ic_\t(\Delta'))^\bullet =\Gamma(\ic_\t(\Delta'))^\bullet$,
the homomorphism
\begin{equation}
\label{inj 3.5}
\H^i(\Gamma(\ic_\t(\Delta(\eta{\prec})))^\bullet) =
\H^i(\Gamma(u_*\ic_\t(\Delta'(\eta{\prec})))^\bullet)\lra
\H^i(\Gamma(u_*\ic_\t(\Delta'))^\bullet)
\end{equation}
is injective for every $i\in\Z$ by Lemma~\ref{Lem 3.4}.

By applying the dualizing functor to the homomorphism
\begin{equation}
\bar\delta_{\Delta'/\Delta}\cl
u_*\ic_\t(\Delta')^\bullet\lra\ic_\t(\Delta)^\bullet\;,
\end{equation}
we get a homomorphism
\begin{equation}
\D(\bar\delta_{\Delta'/\Delta})\cl
\D(\ic_\t(\Delta))^\bullet\lra\D(u_*\ic_\t(\Delta'))^\bullet\;.
\end{equation}
By \cite[Cor.2.12]{Ishida3} and Theorem~\ref{Thm 3.2},
$\D(\ic_\t(\Delta'))^\bullet$ is quasi-isomorphic to
$\ic_\t(\Delta')^\bullet$ in $\CGEM(\Delta')$.
Hence, by Lemma~\ref{Lem 2.3}, (1) and Proposition~\ref{Prop 2.6},
$K^\bullet :=\D(u_*\ic_\t(\Delta'))^\bullet$ is quasi-isomorphic
to $u_*\ic_\t(\Delta')^\bullet$.
By the injectivity of (\ref{inj 3.5}), the homomorphism
\begin{equation}
\label{inj K}
\H^i(\Gamma(K|\Delta(\eta{\prec}))^\bullet)\lra
\H^i(\Gamma(K)^\bullet)
\end{equation}
is injective for every $i\in\Z$.

On the other hand, $L^\bullet :=\D(\ic_\t(\Delta))^\bullet$ is
quasi-isomorphic to $\ic_\t(\Delta)^\bullet$ by
\cite[Cor.2.12]{Ishida3} and Theorem~\ref{Thm 3.2}.
Hence the theorem is equivalent to the injectivity of the
homomorphism
\begin{equation}
\label{inj L}
\H^i(\Gamma(L|\Delta(\eta{\prec}))^\bullet)\lra
\H^i(\Gamma(L)^\bullet)
\end{equation}
for $i\in\Z$.
Let $\Phi :=\bigcup_{\sigma\in\Delta(\eta{\prec})}F(\sigma)$.
Then $\Phi$ is a common subfan of $\Delta$ and $\Delta'$
which contains $\Delta(\eta{\prec})$.
In particular,
$(u_*\ic_\t(\Delta')|\Phi)^\bullet = (\ic_\t(\Delta)|\Phi)^\bullet$.
Hence, by the definition of $\D$, we have
$(L|\Phi)^\bullet = (K|\Phi)^\bullet$ and hence
$(L|\Delta(\eta{\prec}))^\bullet =
(K|\Delta(\eta{\prec}))^\bullet$

We get a commutative diagram
\begin{equation}
\begin{array}{ccc}
\makebox[20pt]{}\H^i(\Gamma(L|\Delta(\eta{\prec}))^\bullet) &
\lra &
 \makebox[10pt]{}\H^i(\Gamma(L)^\bullet) \\
\downarrow &
 &\downarrow \\
\makebox[20pt]{}\H^i(\Gamma(K|\Delta(\eta{\prec}))^\bullet) &
\lra &
\makebox[10pt]{}\H^i(\Gamma(K)^\bullet)
\end{array}
,
\end{equation}
Since $\D(\bar\delta_{\Delta'/\Delta})$ induces the identity
map of $\H^i(\Gamma(L|\Delta(\eta{\prec}))^\bullet)$ and
$\H^i(\Gamma(K|\Delta(\eta{\prec}))^\bullet)$,
the injectivity of (\ref{inj L}) follows from that of
(\ref{inj K}).
\QED

The following theorem is essentially equal to
\cite[Thm.4.2]{Oda2}.
We write here the proof in our notation in order to show
that we need not use the corresponding toric varieties.


\begin{Thm}
\label{Thm 3.6}
Let $\tilde\Delta$ be a complete simplicial fan of
$N_\R$ and $\gamma$ a one-dimensional cone in $\tilde\Delta$.

Let $\Delta$ be the subfan
$\tilde\Delta\setminus\tilde\Delta(\gamma{\prec})$ of
$\tilde\Delta$.
Then
\begin{equation}
\H^p(\Gamma(\ic_\t(\Delta))^\bullet)_q =\{0\}
\end{equation}
for all $p\not= q + r$.
\end{Thm}

\Proof
Since $\Delta(\gamma{\prec})$ is star closed in $\tilde\Delta$,
there exists a short sequence
\begin{equation}
\label{ic short exact}
0\ra\Gamma(\ic_\t(\tilde\Delta(\gamma{\prec})))^\bullet
\lra\Gamma(\ic_\t(\tilde\Delta))^\bullet
\lra\Gamma(\ic_\t(\Delta))^\bullet\ra 0\;.
\end{equation}
We consider the homogeneous degree $q$-part
\begin{equation}
\label{long exact}
\begin{array}{cccccc}
 & & \ra &
\H^{p-1}(\Gamma(\ic_\t(\tilde\Delta))^\bullet)_q &
\ra &
\H^{p-1}(\Gamma(\ic_\t(\Delta))^\bullet)_q \\
\ra &
\H^p(\Gamma(\ic_\t(\tilde\Delta(\gamma{\prec})))^\bullet)_q &
\ra &
\H^p(\Gamma(\ic_\t(\tilde\Delta))^\bullet)_q &
\ra &
\H^p(\Gamma(\ic_\t(\Delta))^\bullet)_q \\
\ra &
\H^{p+1}(\Gamma(\ic_\t(\tilde\Delta(\gamma{\prec})))^\bullet)_q &
\ra
\end{array}
\end{equation}
of the long exact sequence obtained by (\ref{ic short exact})
for each integer $q$.

Since $\tilde\Delta$ is a simplicial complete fan, so is
the fan $\tilde\Delta[\gamma]$ of the $(r-1)$-dimensional
space $N'_\R$.

By the isomorphism
\begin{equation}
h(\tilde\Delta,\gamma)\cl
\ic_\t(\tilde\Delta(\gamma{\prec}))^\bullet\simeq
\epsilon_\gamma(\ic_\t(\tilde\Delta[\gamma]))[-1]^\bullet\;,
\end{equation}
we have an isomorphism
\begin{equation}
\Gamma(\ic_\t(\tilde\Delta(\gamma{\prec})))^\bullet\simeq
\Gamma(\ic_\t(\tilde\Delta[\gamma]))[-1]^\bullet\;,
\end{equation}
of graded $A$-modules.

By Theorem~\ref{Thm 3.3} applied for $\tilde\Delta[\gamma]$,
we have
\begin{equation}
\H^p(\Gamma(\ic_\t(\tilde\Delta[\gamma]))[-1]^\bullet)_q =
\H^{p-1}(\Gamma(\ic_\t(\tilde\Delta[\gamma]))^\bullet)_q =
\{0\}
\end{equation}
for $p - 1\not= q + r - 1$,  i.e., for $p\not= q + r$.
Hence
$\H^p(\Gamma(\ic_\t(\tilde\Delta(\gamma{\prec})))^\bullet)_q =\{0\}$
for $p\not= q + r$.

Since $\tilde\Delta$ is a simplicial complete fan of $N_\R$,
we have
\begin{equation}
\H^p(\Gamma(\ic_\t(\tilde\Delta))^\bullet)_q =\{0\}
\end{equation}
for $p\not= q + r$ by Theorem~\ref{Thm 3.3}.

If $p\not= q + r - 1, q + r$, then we have
$\H^p(\Gamma(\ic_\t(\Delta))^\bullet)_q =\{0\}$ by the
long exact sequence.
We have
$\H^{q+r-1}(\Gamma(\ic_\t(\Delta))^\bullet)_q =\{0\}$
for every $q$, since the homomorphism
\begin{equation}
\label{iso to the shift}
\H^{q+r}(\Gamma(\ic_\t(\tilde\Delta(\gamma{\prec})))^\bullet)\ra
\H^{q+r}(\Gamma(\ic_\t(\tilde\Delta))^\bullet)
\end{equation}
is injective by Theorem~\ref{Thm 3.5}
\QED


\begin{Lem}
\label{Lem 3.7}
Under the same assumption as the above theorem,
the homomorphism
\begin{equation}
\H^{q+r}(\Gamma(\ic_\t(\tilde\Delta))^\bullet)_q\lra
\H^{q+r}(\Gamma(\ic_\t(\Delta))^\bullet)_q
\end{equation}
induced by the natural homomorphism
$\ic_\t(\tilde\Delta)^\bullet\ra\ic_\t(\Delta)^\bullet$
is surjective.
\end{Lem}

\Proof
This lemma follows from the long exact sequence
(\ref{long exact}) in the proof of the theorem, since
$\H^{q+r+1}(\Gamma(\ic_\t(\tilde\Delta(\gamma{\prec})))^\bullet)_q =\{0\}$.
\QED


\section{The diagonal theorems for a complete fan and a cone}
\setcounter{equation}{0}

In this section, we prove two diagonal theorems.


\begin{Thm}[The first diagonal theorem]
		\label{Thm 4.1}
Let $\Delta$ be a complete fan of $N_\R$.
Then
\begin{equation}
\H^p(\Gamma(\ic(\Delta))^\bullet)_q =\{0\}
\end{equation}
for $p, q\in\Z$ with $p\not= q+r$.
\end{Thm}

\Proof
Let $\Sigma\ra\Delta$ be a barycentric subdivision.
Since $\Sigma$ is a simplicial fan, $\ic(\Sigma)^\bullet$
is quasi-isomorphic to $\ic_\t(\Sigma)^\bullet$ by
Theorem~\ref{Thm 3.2}.
Since $\Delta$ is complete, so is $\Sigma$.
Hence $\H^p(\Gamma(\ic(\Sigma))^\bullet)_{q} =\{0\}$ for
$p\not= q + r$ by Theorem~\ref{Thm 3.3}.
We get the theorem, since
\begin{equation}
\dim_\Q\H^p(\Gamma(\ic(\Delta))^\bullet)_{q}\leq
\dim_\Q\H^p(\Gamma(\ic(\Sigma))^\bullet)_{q}
\end{equation}
for any $p, q$ by Theorem~\ref{Thm 2.8}.
\QED


\begin{Thm}
		\label{Thm 4.2}
Let $\Delta$ be a finite fan which may not be complete.
Assume that there exists a complete fan $\tilde\Sigma$ and
a one-dimensional cone $\gamma\in\tilde\Sigma$ such that
$\Sigma :=\tilde\Sigma\setminus\tilde\Sigma(\gamma{\prec})$ is a
barycentric subdivision of $\Delta$.
Then
\begin{equation}
\H^p(\Gamma(\ic(\Delta))^\bullet)_q = 0
\end{equation}
for $p\not= q+r$.
\end{Thm}

\Proof
By Theorem~\ref{Thm 3.6}, we have
$\H^p(\Gamma(\ic(\Sigma))^\bullet)_q = 0$ for $p\not= q+r$.
Then the lemma is a consequence of Theorem~\ref{Thm 2.8}.
\QED


\begin{Thm}[The second diagonal theorem]
		\label{Thm 4.3}
Let $\pi\subset N_\R$ be a cone of dimension $r$.
Then
\begin{equation}
\label{vanishing}
\H^p(\Gamma(\ic(F(\pi)\setmin\{\pi\}))^\bullet)_q =\{0\}
\end{equation}
unless $p + q\geq 0, p\not= q + r - 1$ or
$p + q\leq -1, p\not= q + r$.
\end{Thm}

\Proof
Since $0\leq r_\sigma\leq r-1$ for $\sigma\in F(\pi)\setmin\{\pi\}$,
$\Gamma(\ic(F(\pi)\setmin\{\pi\}))_q^p =\{0\}$
unless $0\leq p\leq r-1$ and $-r\leq q\leq 0$
by \cite[Prop.2.11]{Ishida3}.
In order to prove the theorem, it is sufficient to show
the vanishing (\ref{vanishing}) for $p, q$ with $p\not= q + r - 1$
and $p + q\geq 0$, since then Theorem~\ref{Thm 2.1} implies the
vanishing (\ref{vanishing}) for $p, q$ with $p\not= q + r$
and $p + q\leq -1$.

We denote simply by $\ic(\pi)^\bullet$ the complex
$\ic(F(\pi))(\pi)^\bullet$ in $\CGM(A)$.
Since
\begin{eqnarray}
\ic(\pi)^\bullet & = &
 \grt^{\geq1}(\i_\pi^\circ(\ic(F(\pi)))[-1])^\bullet \\
 & = &
\grt^{\geq1}(\Gamma(\ic(F(\pi)\setmin\{\pi\}))[-1])^\bullet\\
 & = &
(\grt^{\geq0}\Gamma(\ic(F(\pi)\setmin\{\pi\})))[-1]^\bullet\;,
\end{eqnarray}
we have
\begin{equation}
\H^p(\Gamma(\ic(F(\pi)\setmin\{\pi\}))^\bullet)_q =
\H^{p+1}(\ic(\pi)^\bullet)_q
\end{equation}
for $p + q\geq 0$.
Hence it is sufficient to show that
$\H^p(\ic(\pi)^\bullet)_q = 0$ for $p\not= q + r$.

We take a one-dimensional cone $\gamma$ of $N_\R$ such that
$-\gamma$ intersects the interior of $\pi$.
We set
\begin{equation}
\Delta := (F(\pi)\setminus\{\pi\})\cup
\{\gamma+\sigma\scl\sigma\in F(\pi)\setminus\{\pi\}\}
\end{equation}
and
\begin{equation}
\tilde\Delta :=\Delta\cup\{\pi\}\;.
\end{equation}
Then $\tilde\Delta$ is a complete fan.

Let $\tilde\Sigma$ be a barycentric subdivision of
$\tilde\Delta$ and let $\eta\in\tilde\Sigma$ be the
one-dimensional cone which intersects the interior of $\pi$.
Then $\Sigma :=\tilde\Sigma\setminus\tilde\Sigma(\eta{\prec})$
is a barycentric subdivision of $\Delta$.

By Theorems~\ref{Thm 4.1} and \ref{Thm 4.2},
$\H^p(\Gamma(\ic(\tilde\Delta))^\bullet)_q$ and
$\H^p(\Gamma(\ic(\Delta))^\bullet)_q$ are zero unless
$p = q + r$.

By the exact sequence
\begin{equation}
0\ra\ic^\bullet(\pi)\ra\Gamma(\ic(\tilde\Delta))^\bullet
\ra\Gamma(\ic(\Delta))^\bullet\ra 0\;
\end{equation}
we get the exact sequence
\begin{equation}
\begin{array}{ccccccc}
0 & \ra & \H^{q+r}(\ic^\bullet(\pi))_q & \ra &
\H^{q+r}(\Gamma(\ic(\tilde\Delta))^\bullet)_q &
\mathop{\ra}\limits^{\varphi} & \\
 & & \H^{q+r}(\Gamma(\ic(\Delta))^\bullet)_q & \ra &
\H^{q+r+1}(\ic^\bullet(\pi))_q & \ra & 0
\end{array}
\;,
\end{equation}
while we have $\H^p(\ic^\bullet(\pi))_q =\{0\}$
for $p\not= q + r, q + r + 1$.
Hence it is sufficient to show the surjectivity of $\varphi$
for each $q\in\Z$.

Then we get a commutative diagram of canonical homomorphisms
\begin{equation}
\begin{array}{ccc}
\H^{q+r}(\Gamma(\ic(\tilde\Sigma))^\bullet)_q &
\mathop{\lra}\limits^{\varphi'} &
\H^{q+r}(\Gamma(\ic(\Sigma))^\bullet)_q \\
\hbox{  }\downarrow{\scriptstyle\psi'} & &\hbox{  }\downarrow{\scriptstyle\psi}
\\
\H^{q+r}(\Gamma(\ic(\tilde\Delta))^\bullet)_q &
\mathop{\lra}\limits^{\varphi} & \H^{q+r}(\Gamma(\ic(\Delta))^\bullet)_q
\end{array}
\;.
\end{equation}

Since $\tilde\Sigma$ is simplicial, $\ic(\tilde\Sigma)^\bullet$
is quasi-isomorphic to $\ic_\t(\tilde\Sigma)^\bullet$.
Hence $\varphi'$ in the diagram is surjective
by Lemma~\ref{Lem 3.7}.
Since $\Sigma$ is a barycentric subdivision of $\Delta$, $\psi$
is also surjcetive by Theorem~\ref{Thm 2.8}
Hence $\varphi$ is surjective for every $q$.
\QED

Note that the theorem says that the cohomologies may not
zero only for $(p, q) = (i-1, i-r)$ for $r/2 < i\leq r$
or $(p, q) = (i, i-r)$ for $0\leq i < r/2$.

Here we comment the relation with the results of \cite{Oda3}.

If the cone $\pi$ is of dimension $r$ and $F(\pi)\setminus\{\pi\}$
is a simplicial fan, then $\ic(F(\pi)\setmin\{\pi\})^\bullet$
is quasi-isomorphic to $\ic_\t(F(\pi)\setmin\{\pi\})^\bullet$
by Theorem~\ref{Thm 3.2}.
Hence the second diagonal theorem implies that
\begin{equation}
\label{II for simp}
\H^p(\Gamma(\ic_\t(F(\pi)\setmin\{\pi\}))^\bullet)_q =\{0\}
\end{equation}
unless $p + q\geq 0, p\not= q + r - 1$ or
$p + q\leq -1, p\not= q + r$.

Let $\Pi$ be a simplicial complete fan of $N_\R$.
Suppose a continuous map $h\cl N_\R\ra\R$ satisfies the
following conditions.

(1) $h$ is linear on each cone $\sigma\in\Pi$.

(2) $h$ has rational values on $N_\Q$.

(3) $h$ is strictly convex with respect to the fan $\Pi$,
i.e., $h(x) + h(y)\geq h(x + y)$ and the equality holds
only if $x, y$ are in a common cone of $\Pi$.

We set $N' := N\oplus\Z$ and
$\tilde\sigma :=\{(x,h(x))\scl s\in\sigma\}\subset N'_\R$ for
each $\sigma\in\Pi$, and define
\begin{equation}
\tilde\Pi :=\{\tilde\sigma\scl\sigma\in\Pi\}\;.
\end{equation}
Then $\tilde\Pi$ is a simplicial fan of the $(r+1)$-dimensional
space $N'_\R$ and $\tilde\Pi = F(\pi)\setminus\{\pi\}$ for
the $(r+1)$-dimensional cone
$\pi := \{(x, y)\scl x\in N_\R, y\geq h(x)\}\subset N'_\R$.

We see that the complex $C^\bullet(\tilde\Pi,\tilde\calG_\ell)$ defined
in \cite{Oda3} is isomorphic to $\ic_\t(\tilde\Pi)_{\ell-r-1}^\bullet$
for each $0\leq\ell\leq r + 1$.
Hence by the second diagonal theorem, we have
\begin{equation}
\H^p(\tilde\Pi,\tilde\calG_\ell) :=
\H^p(C^\bullet(\tilde\Pi,\tilde\calG_\ell)) =\{0\}
\end{equation}
except for $(p,\ell) = (i-1, i)$ for $(r+1)/2 < i\leq r + 1$ or
$(p,\ell) = (i, i)$ for $0\leq i < (r+1)/2$.
In particular, it is zero for $(p,\ell) = (i, i)$ with
$r/2 < i$ and $(p,\ell) = (i-1, i)$ with $i < r/2 + 1$
since $r$ and $i$ are integers.

Hence we get the equivalent conditions of \cite[Cor.4.5]{Oda3}
for a (rational) fan $\Pi$ and a rational strictly convex
function $h$.
Note that these conditions for the rational case is a consequence
of the strong Lefschetz theorem (cf.\cite{Oda3}).

However irrrational fans are also treated in \cite{Oda3},
our theory can not be applied for irrational case.


\begin{Cor}
		\label{Cor 4.4}
Let $\rho$ be a nontrivial cone in $N_\R$.
If $i + j\geq 0$ and $i\not= j + r_\rho - 1$ or if
$i + j\leq -1$ and $i\not= j + r_\rho$, then
$\H^i(\i_\rho^*(\ic(F(\rho)\setmin\{\rho\}))^\bullet)_j =\{0\}$.
\end{Cor}

\Proof
The cone $\rho$ is of maximal dimension in the real space
$N(\rho)_\R$ of dimension $r_\rho$.
Since the functor $\i_\rho^*$ for the fan $F(\rho)$ is equal
to $\Gamma$ of this space, this is a consequence of
Theorem~\ref{Thm 4.3}.
\QED


\begin{Cor}
		\label{Cor 4.5}
Let $\rho$ be in a finite fan $\Delta$.
Then we have $\H^i(\ic(\Delta)(\rho)^\bullet)_j =\{0\}$
unless $i + j\geq 1,\; i = j + r_\rho$, while
$\H^i(\i_\rho^*(\ic(\Delta))^\bullet)_j =\{0\}$ unless
$i + j\leq -1,\; i = j + r_\rho$.
\end{Cor}

\Proof
By the construction of $\ic(\Delta)^\bullet$ in
\cite[Thm.2.9]{Ishida3},
$\ic(\Delta)(\rho)^\bullet$ is isomorphic to
$\grt^{\geq 1}(\i_\rho^*(\ic(F(\rho)\setmin\{\rho\}))[-1])^\bullet$
in $\CGM(A(\rho))$, and $\i_\rho^*(\ic(\Delta))^\bullet$ is
quasi-isomorphic to
$\grt_{\leq -1}(\i_\rho^*(\ic(F(\rho)\setmin\{\rho\})))^\bullet$.
Hence this is a consequence of Corollary~\ref{Cor 4.4}.
\QED

\end{document}